\documentclass{article}
\usepackage[utf8]{inputenc}

\usepackage{amsmath}
\usepackage{amssymb}
\usepackage{amsthm}
\usepackage{mathtools}
\usepackage{appendix}

\usepackage{graphicx}
\usepackage{psfrag}
\usepackage{authblk}
\usepackage{color}
\usepackage{caption}
\usepackage{subcaption}
\usepackage{url}

\usepackage{xtab,booktabs}

\newcommand{\N}{\mathbb N}
\newcommand{\R}{\mathbb R}

\newcommand{\Z}{\mathbb Z}

\usepackage{xcolor}
\usepackage{sqrcaps}

\numberwithin{theo}{section}

\newtheorem{lem}[equation]{Lemma}

\usepackage[symbol]{footmisc}

\makeatletter
\newcommand\footnoteref[2]{\protected@xdef\@thefnmark{\ref{#1}}\@footnotemark}
\makeatother

\begin{document}

\begin{center}
\section*{How neural networks learn to classify chaotic time series}

\large{Alessandro Corbetta\footnote[2]{\label{note1}Joint first author}\\}
\small{Department of Applied Physics, \\
Eindhoven University of Technology, Eindhoven, The Netherlands. \\
{\tt a.corbetta@tue.nl}}
\\[2mm]

\large{Thomas Geert de Jong\footnotemark[2]\\}
\small{Faculty of Mathematics and Physics, 
Institute of Science and Engineering, \\
Kanazawa University,
Kanazawa, Japan.\\
{\tt tgdejong@se.kanazawa-u.ac.jp}}
\\[2mm]
\end{center}

\begin{small} \noindent \textbf{Abstract.} Neural networks are
  increasingly employed to model, analyze and control
  non-linear dynamical systems ranging from physics to biology.  Owing to their
  universal approximation capabilities, they regularly
  outperform state-of-the-art model-driven methods in terms of accuracy,
  computational speed, and/or control capabilities. On the other hand,
  neural networks are very often they are taken as black boxes whose
  explainability is challenged, among others, by huge amounts of
  trainable parameters.

  In this paper, we tackle the outstanding issue of analyzing the
  inner workings of neural networks trained to classify 
  regular-versus-chaotic time series. This setting, well-studied in dynamical systems, enables
  thorough formal analyses. We focus specifically on a family of networks dubbed Large Kernel
  Convolutional Neural Networks (LKCNN), recently introduced by
  Boull\'{e} et al. (2021). These non-recursive networks have been
  shown to outperform other established architectures (e.g. residual
  networks, shallow neural networks and fully convolutional networks)
  at this classification task. Furthermore, they outperform ``manual''
  classification approaches based on direct reconstruction of the
  Lyapunov exponent. We find that LKCNNs use qualitative properties of
  the input sequence. In particular, we show that the relation between
  input periodicity and activation periodicity is key for the
  performance of LKCNN models. Low performing models show, in fact,
  analogous periodic activations to random untrained models. This
  could give very general criteria for identifying, a priori, trained
  models that have poor accuracy.

\vspace{3mm}

\noindent {\bf Keywords:} Dynamical systems, chaos, deep learning, convolutional networks, time series, classification, Savitsky-Golay.\\  
\end{small}


\section{Introduction \label{sec:intro}}

During the last decade there has been a strong acceleration in the adoption of machine learning, typically through artificial neural networks, to model, analyze, and control a broad spectrum of dynamical systems~\cite{brunton2022data}. 
Neural networks are nonlinear parametric models satisfying universal approximation principles~\cite{uniapprox}. Typically, they are trained in data-driven contexts to fit pre-annotated data. Beside boosting well-known super-human performance in contexts as automated vision~\cite{wieting2016ICLR}, natural-language processing~\cite{sebe2005machine}, and long-term strategy games~\cite{silver2016mastering}, they are gaining substantial momentum in fundamental research in connection with non-linear/chaotic dynamics. In  these contexts they are regularly surpassing ``traditional'' state-of-the-art model-driven approaches~\cite{tang2020introduction}. \\
%
%

Astounding are the achievements connected to three-body problem long-term forecasts~\cite{breen2020newton}, weather modeling~\cite{kashinath2021physics},  fluid turbulence closures~\cite{duraisamy2019turbulence,ortali2022numerical},  measurements~\cite{corbetta2021deep}, and control~\cite{rabault2019artificial}. Machine learning methodologies allowed progress in scientific context-agnostic issues, e.g., model-free chaotic dynamics predictions~\cite{pathak2018model}, as well as the problem dealt with in this paper: the classification of time series between chaotic and regular. Dynamics that can bifurcate from regular to chaotic and vice versa are present in every scientific discipline: astronomy, biology, meteorology, physics etc.~\cite{peitgen2004chaos,strogatz2018nonlinear,celletti2010stability,lorenz1963deterministic, kuznetsov1998elements,guckenheimer2013nonlinear,broer2011dynamical}. While the issue of their classification is thus longstanding (e.g.~\cite{kuo1992prediction,lai2003recent,mukherjee1997nonlinear, panday2021machine, principe1992prediction,woolley2010modeling, boulle2020classification}), recent works by Boull\'{e} et al.~\cite{boulle2020classification} showed that convolutional neural networks, specifically with large kernels, can be strikingly successful at the task.  In general, convolutional networks~\cite{goodfellow2016deep} have proven suitable for time series classification~\cite{wang2017time} and inference~\cite{corbetta2021deep}. Building on this approach, Boull\'{e} et al. proposed to apply Large Kernel Convolutional Neural Networks, LKCNN. Large kernel means that the kernel size of the convolutional layers is large in relation to the input sequence length. Boull\'{e} et al. show that LKCNNs boast highest performance in comparison with other established machine learning approaches such as residual, fully convolutional, multilayer perceptrons and shallow networks~\cite{wang_github,wang2017time}. Besides they manage generalization properties to data outside the training set. 

This creates the case for the two-fold analysis of this paper, only possible because of two factors: the relative technical simplicity of LKCNNs and the vast understanding on chaotic maps. Thus, this paper contributes thorough comparisons in terms of accuracy with traditional approaches, currently missing and a formal analysis of the inner workings of LKCNNs that allow such performance. 

In general, neural networks allow no straightforward insight into their internal mechanics, especially when it comes to large and complex networks~\cite{bengio2021deep}. Thus, while networks are ubiquitously employed as (effective) data-driven black boxes, understanding the inner workings, e.g. identifying the features that they leverage on, could become a crucial component towards new discoveries. \\

In this paper, we analyze how neural networks successfully manage the longstanding challenge of classifying discrete-time series produced by dynamical systems, identifying regular vs. chaotic motions. Here, we will take the approach of identifying chaos through sensitive dependence on initial conditions by computing the Lyapunov exponent which measures the average rate of exponential departures of small perturbations~\cite{benettin1980lyapunov, parker2012practical}. \newpage

Formally, we consider the following classification problem:  Let
$f: \R \rightarrow \R $ be a smooth map recursively generating the bounded time
series $\{ x_{1},x_{2},\ldots\}$ as
\begin{align}
x_{n+1} = f(x_n), \label{eq:f}
\end{align}
for a given initial condition $x_0$. We aim at classifying finite-length time sequences between \textit{chaotic} or \textit{regular} without knowledge of the analytical expression of $f$. In other terms, given finite sequences 
\begin{equation*}
x^\ell = ( x_{1},x_{1},\dots, x_{N}) \in \R^N,
\end{equation*}
we target accurate classifiers $G$, such that
\begin{gather}
\begin{aligned}
G: \R^N &\rightarrow L := \{ {\rm chaotic , \;  regular }\},\\
x^\ell &\mapsto \ell \in L.
\end{aligned} \label{main:problem}
\end{gather}
Analogously to~\cite{boulle2020classification}, we shall consider $N=500$. Note that if the function $f$ and its derivative are known analytically, the classification labels in $L$ can be immediately determined by estimating the Lyapunov exponents of $f$. More specifically, the Lyapunov exponent, $\lambda$, corresponding to  $ \{ x_{n} \}$ generated by $f$ is defined as
\begin{align}
    \lambda = \lim_{n \rightarrow \infty} \frac{1}{n}\sum^{n-1}_{i=0} \log(|f'(x_i)|), \label{def:lya}
\end{align}
if the limit exists. Bounded sequences are chaotic when $\lambda >0$, as attractors with $\lambda>0$ can only be bounded if a type of folding process merges widely separated trajectories (e.g.~\cite{parker2012practical,wolf1985determining}). Conversely, when only data in the form of sequences ($x^\ell$) is available, direct approaches which reconstruct the Lyapunov exponent by estimating $f'$ could be used. These approaches, however, require a sufficiently large $N$ (e.g.~\cite{vastano1986comparison,bryant1990lyapunov}). Our experiments show that our choice $N=500$ is large enough to capture the dynamics and short enough that direct approaches are error prone. \\

Here, we analyze the inner workings of LKCNNs when classifying time series from well-known 1D maps: the logistic map and the sine-circle map. LKCNNs are conceptually simple: this enables us to perform in-depth experimental and formal analyses. Additionally, we compare LKCNNs performance with an estimate of $\lambda$ via a local polynomial reconstruction of $f'$ using Savisky-Golay filtering~\cite{savitzky1964smoothing}. We restrict to a data set made of non-periodic trajectories with Lyapunov exponents close to zero as these are the most challenging to classify:  small perturbations can lead in fact to misclassification. For these, the derivative along the orbits typically oscillate around 1. Overall, we show that:
    \begin{itemize}
\item[-] LKCNNs use qualitative properties of the input sequence as they outperform direct reconstruction of the Lyapunov exponent from the time series. Our experiment is performed on a restricted data set with only non-periodic sequences. These are  hard to classify as the average rate of separation between consecutive time-steps is almost constant.
\item[-] LKCNNs map periodic inputs to periodic activation with exception of the dense layers. We can capture this mapping in a two dimensional matrix which we refer to as the \textit{period matrix}. The period matrix is model-independent over generic untrained LKCNNs and corresponds to a property of the architecture. We prove this property in a limit condition.  
\item[-] Grouping trained LKCNNs by period matrices, we observe a single period matrix which correlates with low performance models. This period matrix is equal to the period matrix of generic untrained models. 
    \end{itemize}
This last insight might also be useful for addressing more general models and settings.

The remainder of this paper is organized as follows. In Section~\ref{sec:dynam}, we
present the dynamical systems for the classification problem, datatsets and performance metrics. In Section~\ref{sec:SG}, we present a method to directly reconstruct the Lyapunov exponent from the times series with Savitsky-Golay polynomial filtering. In
Section~\ref{sec:LKCNN}, we revise LKCNN, present the architecture,
explain the choice of hyperparameters and data sets.  In
Section~\ref{sec:res}, we present the main results related to Lyapunov
reconstruction with Savitsky-Golay and periodic activation. A final
discussion closes the paper.

\section{Time series, data sets, performance metrics \label{sec:dynam}}

\subsection{Time series: the logistic and sine circle map}
\begin{figure}[ht]
     \centering
     \begin{subfigure}{.45\textwidth}
         \centering
         \includegraphics[width=\linewidth]{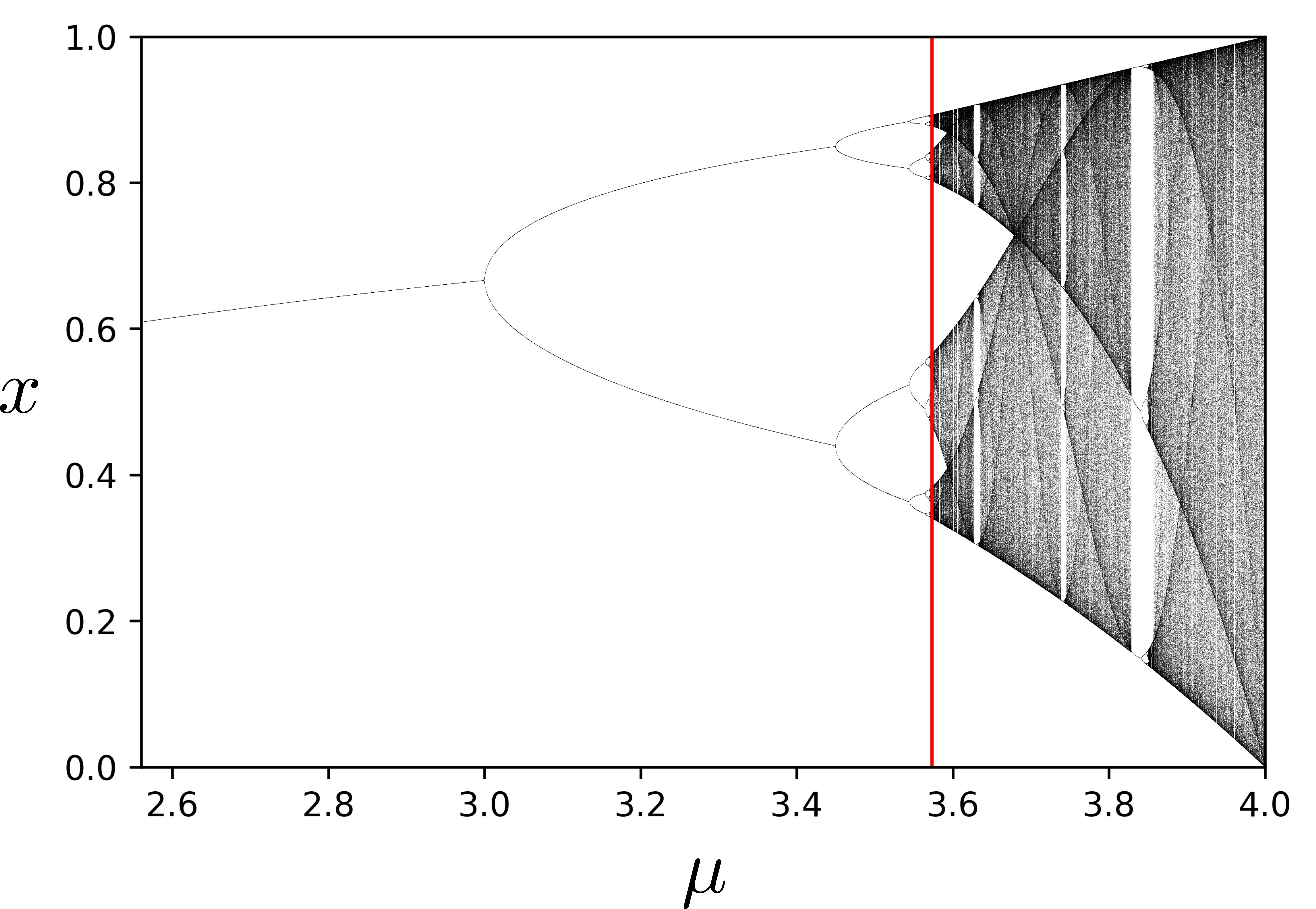}
         \caption{Logistic map \label{fig:log_bifu}}
     \end{subfigure}
     \begin{subfigure}{.45\textwidth}
         \centering
         \includegraphics[width=\linewidth]{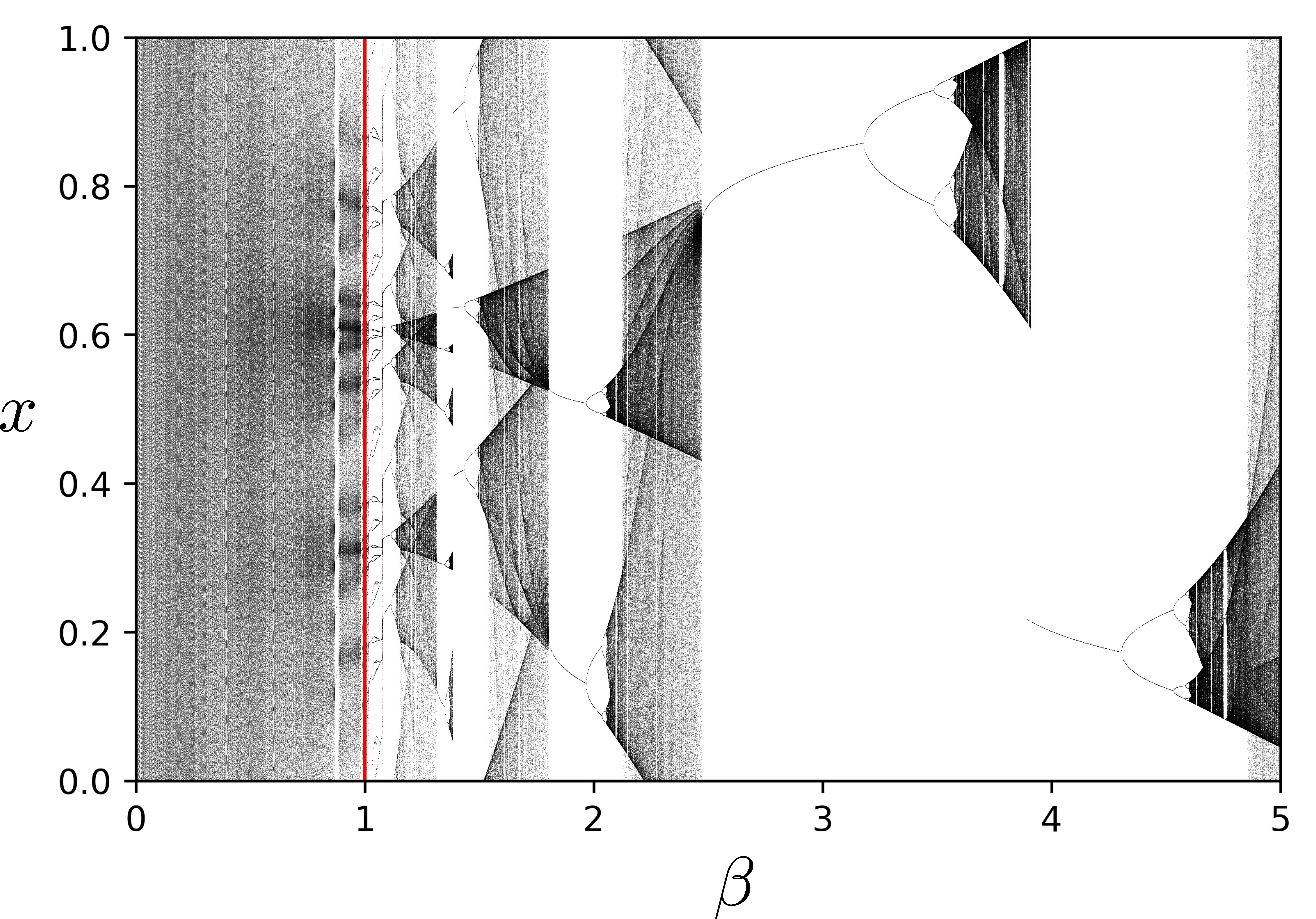}
         \caption{Sine-circle map \label{fig:sine_bifu}}
     \end{subfigure}
     \caption{Bifurcation diagrams: In (a) the parameter corresponding to the onset of chaos is indicated by a red line. In (b) the parameter corresponding to the transition from invertible to non-invertible is indicated by the red line. Invertible maps do not exhibit chaotic dynamics.  \label{fig:bifus}}
\end{figure}

As in~\cite{boulle2020classification}, we analyze LKCNN considering data sets built using two well-known maps: the logistic map and the sine-circle map. We briefly review some crucial features of the time series that these maps generate. The logistic map~\cite{robert1976simple} is given by

\begin{align}
  f(x) = \mu x(1-x)  , \; \mu \in [0,4], \; x_0  = 0.5, \label{def:logmap}
\end{align}
where $\mu$ is the bifurcation parameter. Observe that $x_n \in [0,1]$ for all $\mu$. The bifurcation diagram corresponding to long-term evolution of the orbits is given in Figure~\ref{fig:log_bifu}. The bifurcation diagram exhibits orbits doubling in period with geometric rate. After the period-doubling accumulation point, $\mu \approx 3.56995$, chaos onsets and the map exhibits alternation of regular and chaotic behavior~\cite{feigenbaum1978quantitative}.

The second map considered is the sine-circle map \cite{hilborn2000chaos,boyland1986bifurcations,herman1979conjugaison}:
\begin{align}
f(\theta) = \theta + \Omega - \frac{\beta}{2\pi} \sin{2 \pi \theta} \mod{1}, \; \;  \theta_0=0.5\;. \label{def:sinemap}  
\end{align}
The sine-circle map also exhibits a transition from regular to chaotic dynamics as the bifurcation parameter $\beta$ changes. For $\beta<1$ only periodic and quasi-periodic trajectories occur since the sine-circle map is invertible. This implies that the folding dynamics which characterizes chaotic dynamics cannot occur~\cite{hilborn2000chaos}. As in~\cite{boulle2020classification} we consider $\Omega = 0.606661$, $\beta \in [0,5]$. For this particular $\Omega$ value, chaotic dynamics occur right after $\beta=1$~\cite{hilborn2000chaos}. The bifurcation diagram corresponding to long-term evolution of the orbits is given in Figure \ref{fig:sine_bifu}.

\subsection{Data sets  \label{sec:data}}

\begin{figure}[ht]
     \centering
     \begin{subfigure}{.45\textwidth}
         \centering
         \includegraphics[width=\linewidth]{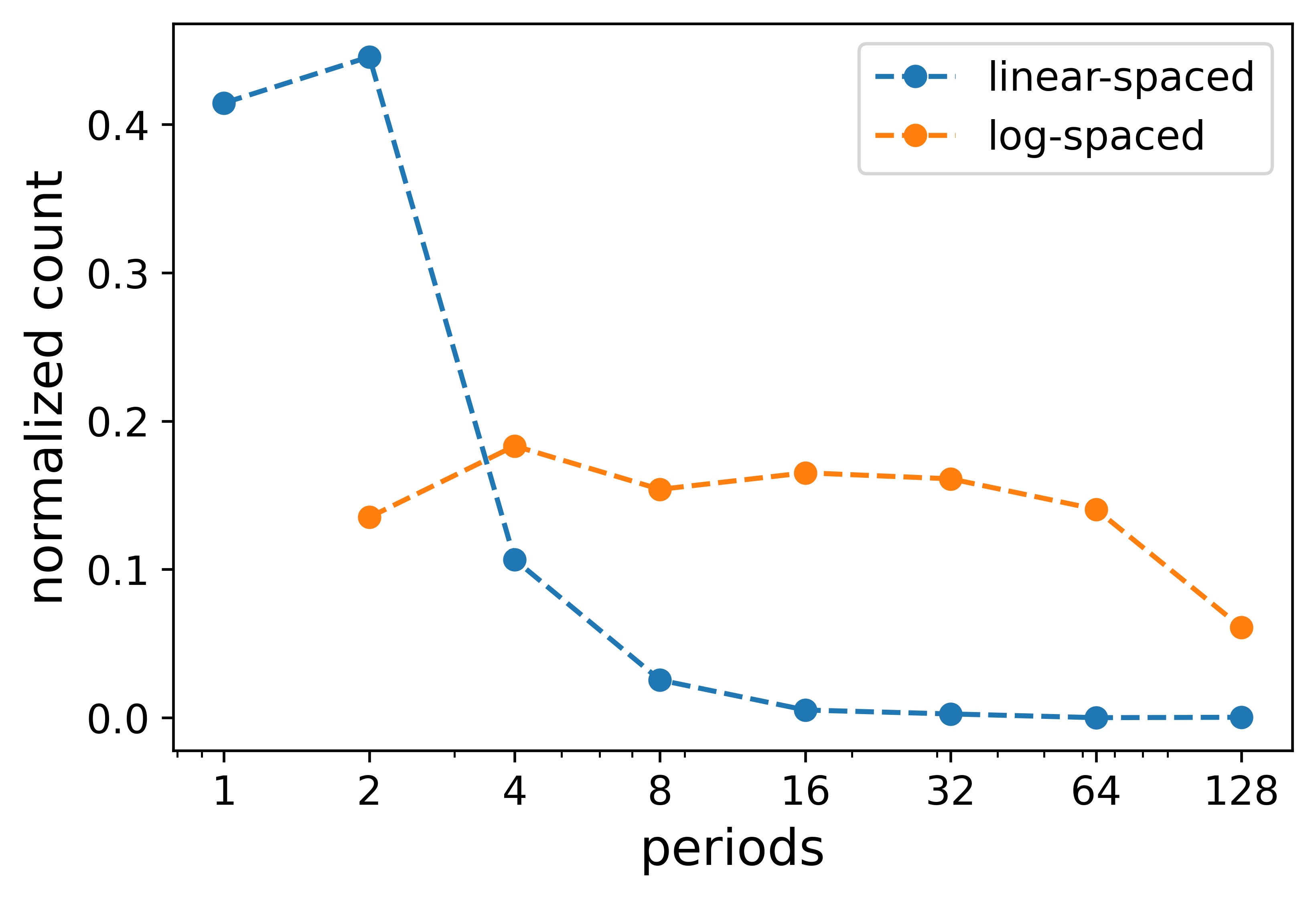}
         \caption{Before onset of chaos}
     \end{subfigure}
     \begin{subfigure}{.45\textwidth}
         \centering
         \includegraphics[width=\linewidth]{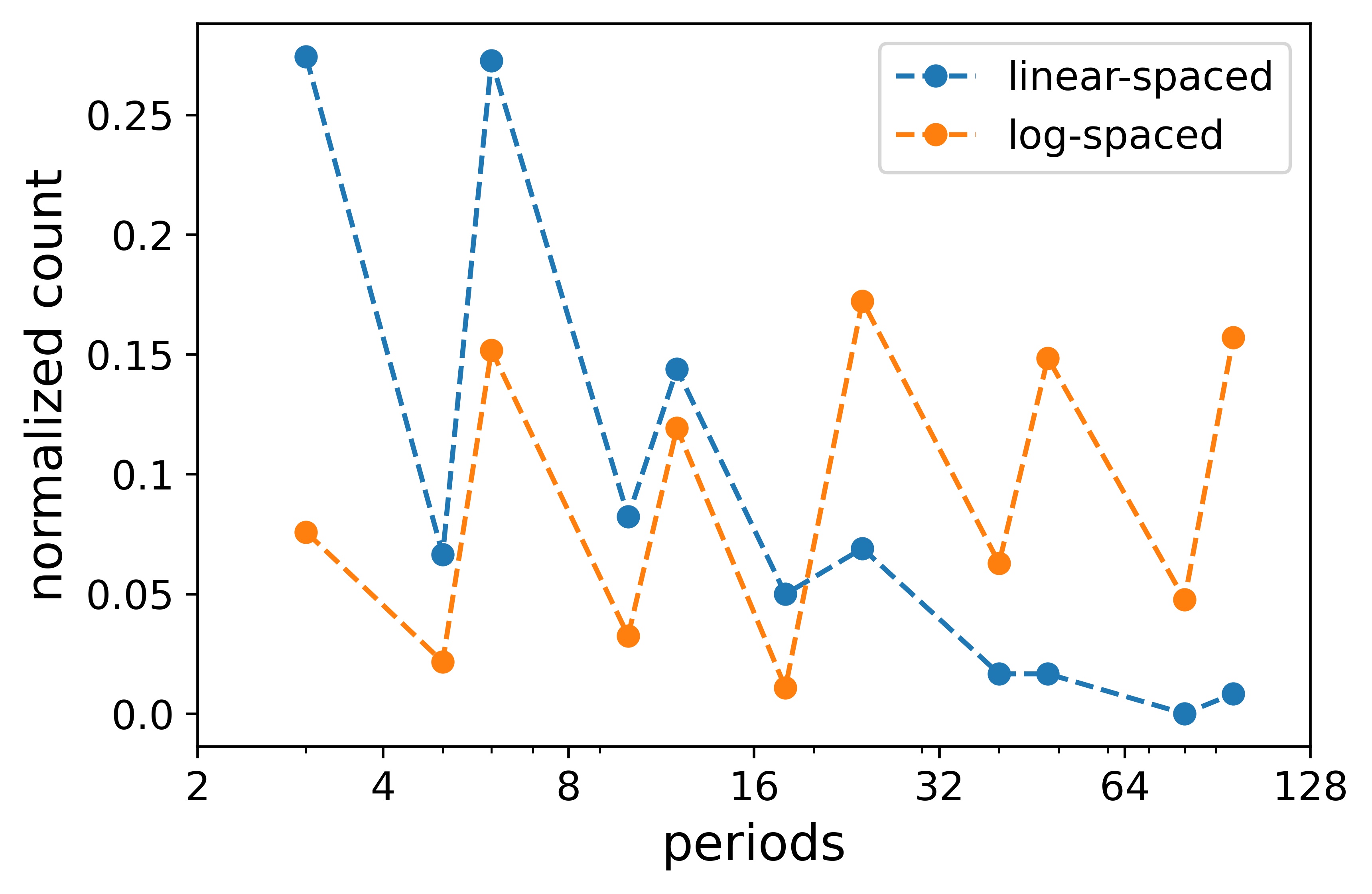}
         \caption{After onset of chaos}
     \end{subfigure}
     \caption{Logistic map period distributions for linear- and log-spaced data: We consider periodic orbits before and after the onset of chaos. The period distribution for the log-spaced data set is more uniform than the linear-spaced data set. \label{fig:period_distribution} }
\end{figure}

Given a discrete parameter set we consider trajectories corresponding to long-term evolution of the time series generated by \eqref{def:logmap} and \eqref{def:sinemap}. We use log-spaced data for the logistic data set as chaotic and regular trajectories are better balanced as a consequence of the geometric doubling over the parameter space. In~\cite{boulle2020classification}, the linear-spaced data set is considered for the logistic map. This data set is biased towards regular trajectories, see Table~\ref{tab:data_stat}. We observe that periods of periodic orbits are much more uniformly distributed in the log-spaced data set than in the linear-spaced data set, Figure~\ref{fig:period_distribution}. Moreover, periodic orbits with low orbit periods are over-represented in the linear data set. 

For the sine-circle map we consider a linear-spaced parameter set.  The regular trajectories can exhibit periodic and non-periodic behavior. The chaotic trajectories are all non-periodic. We note that the sign of the short time Lyapunov exponent is a good predictor of the sign of the converged Lyapunov exponent if the absolute Lyapunov exponent is sufficiently large. The sine-circle data set is well-balanced between chaotic and regular trajectories, see Table \ref{tab:data_stat}. We can heuristically explain why regular orbits make up at least $\approx$ 60\%. For $\beta \in (0,1)$ the motion is periodic or quasi-periodic \cite{hilborn2000chaos} and from Figure \ref{fig:sine_bifu} we observe that period 1, period 2 and period 4 cover a parameter range close to length 2.

\begingroup
\begin{table}
\centering
\begin{tabular}{lrrr}
\toprule
Data sets & Chaotic & \multicolumn{2}{c}{Regular} \\
 \cmidrule(lr){3-4}
&  & Periodic  &  Non-periodic \\
\midrule
Logistic, linear-spaced &  27\%  & 72 \% & 1 \% \\
Logistic, log-spaced & 46\% &  45 \% & 9 \%  \\
Sine-circle, linear-spaced &  41\% & 49\% & 10 \% \\
\bottomrule
\end{tabular}
\caption{Subdivision data sets: The logistic, log-spaced data set is a more balanced data set in comparison to the logistic linear-spaced data set.  For the log-spaced data the parameters are logarithmically spaced at the onset of chaos. We note that regular non-periodic orbits for the logistic map typically correspond to orbits which accumulate around orbits of low period whereas regular non-periodic orbits for the sine-circle map can exhibit quasi-periodic dynamics which covers a large domain of the phase space. \label{tab:data_stat}}
\end{table}
\endgroup

\subsection{Performance metrics}
The performance is assessed using the classification accuracy which is defined as
\[
{\rm Accuracy} := \frac{T_{\rm regular} + T_{\rm chaotic}}{T_{\rm regular} + T_{\rm chaotic} + F_{\rm regular} + F_{\rm chaotic}},
\]
where $T_{\rm regular} = $ true regular predictions, $T_{\rm chaotic}=$ true chaotic predictions, $F_{\rm regular}=$ false regular predictions and $F_{\rm chaotic} =$ false chaotic predictions. Our test sets are well-balanced hence we do not need to consider alternatives \cite{olson2008advanced, brodersen2010balanced}. We will also consider a precision measure for chaotic and regular orbits:
\begin{equation}
 P_{\rm chaotic} := \frac{T_{\rm chaotic}}{T_{\rm chaotic} + F_{\rm chaotic}}, \qquad {P}_{\rm regular} := \frac{T_{\rm regular}}{T_{\rm regular} + F_{\rm regular}}.   \label{eq:precision}
\end{equation}

\section{Direct reconstruction of Lyapunov exponents via Savitsky-Golay polynomial filtering \label{sec:SG}} 

We can reconstruct the Lyapunov exponent directly from the times series by approximating derivatives. Here, we consider Savitsky-Golay as the fluctuations occurring for Lyapunov exponent close to zero are smoothed out which makes it outperform the classical method~\cite{wolf1985determining,rosenstein1993practical}.

Savitsky-Golay is used to obtain the graph of $f$ \eqref{eq:f}. The coefficients of the Savitsky-Golay approximation can be used to compute the derivative of $f$. Standard Savitsky-Golay works on an equidistant grid~\cite{savitzky1964smoothing} but the extension to a non-equidistant grid is straightforward. We will briefly outline the method.\\
The approximation of the graph of $f$ is described by $(x_i,x_{i+1})$. We sort these points with respect to $x_i$. Denote the reordered points by $(\overline{x}_i, y_i)$. Consider $k$ consecutive points $(\overline{x}_i,y_i)$ through which we interpolate an  $m$-degree polynomial with coefficients $\mathbf{c} \in \R^m$. For convenience we will consider $i=1, \cdots ,k$. Take $\kappa = (k-1)/2$. Define the non-equidistant design matrix by $A_{ij} = (\overline{x}_i-\overline{x}_{\kappa})^j$. Define $\mathbf{y}=(y_1, \cdots, y_{k} )^T$. The set of equations that need to be solved are 
\[
A \mathbf{c} = \mathbf{y}.
\] 
From $\mathbf{c}$ we can obtain the derivatives at $x_i$ by straightforward calculus. Finally, we log-average over these derivatives to obtain an estimate for the Lyapunov exponent.

The graph of the attractor will generally consist of components which are disconnected from each other. So we perform Savitsky-Golay polynomial filtering on each component separately. Finally, we note that this procedure does not directly extend to higher dimensions.

\section{Large Kernel Convolutional Neural Networks (LKCNN) \label{sec:LKCNN}}
We review here the core features of Large Kernel Convolutional Neural Networks (LKCNN) as introduced in~\cite{boulle2020classification}. A Large Kernel Convolutional Neural Network  has convolutional layers with  kernel size large in relation to its feature channels. Conversely, standard convolutional networks usually consider a large number of feature channels in relation to the kernel size~\cite{ismail2019deep, wang2017time}. In this paper, as in \cite{boulle2020classification}, we use LKCNN for binary classification of sequences. 

\subsection{Architecture}

\begin{figure}[ht]
    \centering
    \includegraphics[width=10cm]{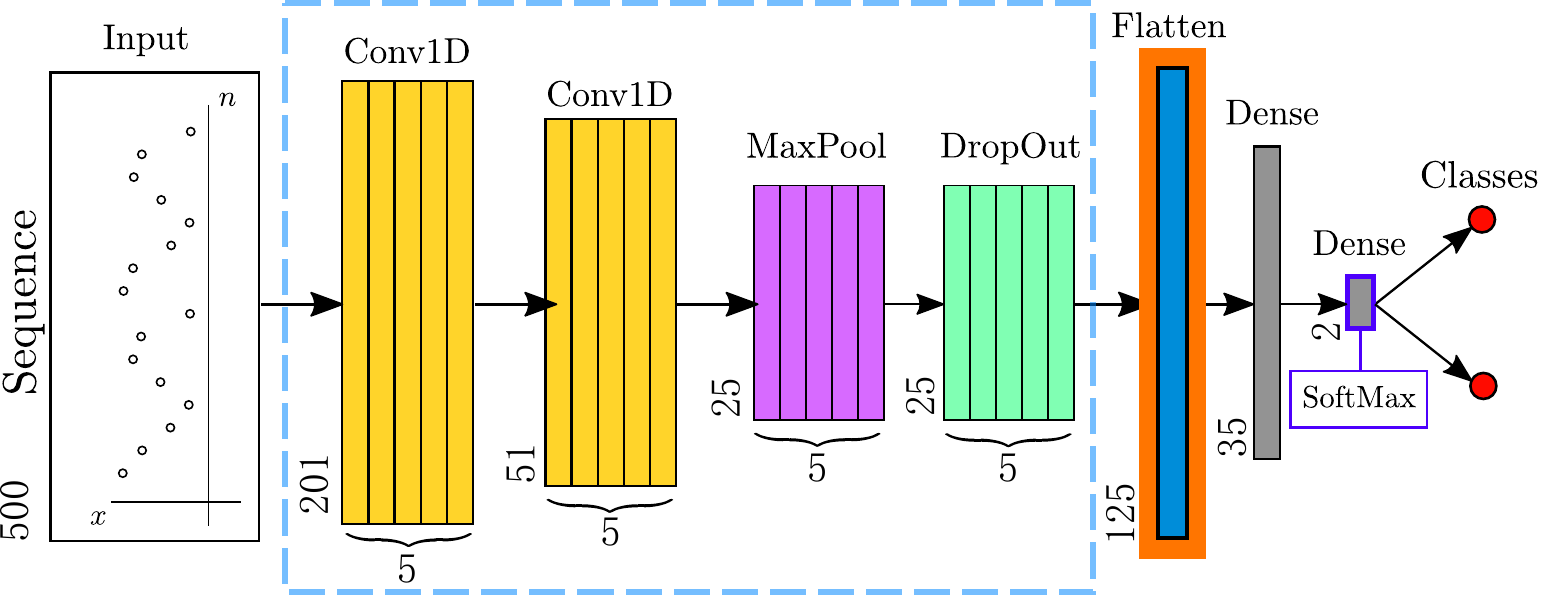}
    \caption{LKCNN architecture: The input time series $x$ is classified as chaotic or regular. LKCNN have a typical CNN architecture but the convolutional layers have a large kernel size. The two Conv1D each have 5 filters with kernel size 100. The loss function used is cross-entropy. The internal analysis of the network will focus on the effect of the input on the activation at the flatten layer, highlighted in orange. Specifically, we will see how the activation of the layers in the blue dashed rectangle is mapped to an activation at the flatten layer. }
    \label{fig:LKCNN_arch}
\end{figure}

We report in Figure~\ref{fig:LKCNN_arch} the LKCNN architecture we employ to tackle~\eqref{main:problem}. 
Specifically, we stack the following layers in a feed-forward fashion:
\begin{itemize}
\item two convolutional layers (kernel size: $100$, large in relation to the signal length; stride 2; relu activation);
\item maxpooling layer (pool size 2);
\item dropout regularization layer (rate 0.5); 
\item dense layer (sigmoid activation);
\item final dense layer with softmax activation over the two classes in $L$. In other words, the network outputs a probability vector $(p_{\rm regular}, p_{\rm chaotic}) = (p_{\rm regular}, 1-p_{\rm regular})$.
\end{itemize}
We consider an input sequence as chaotic if $p_{\rm chaotic} > 0.5$ and regular otherwise.  We employ cross-entropy as training loss. 

The overall architecture is the same as in~\cite{boulle2020classification} but with different hyperparameters. In fact, we use exclusively the Lyapunov exponent~\eqref{def:lya} to assign the classification labels on the training set. In particular, we estimate the Lyapunov exponent employing orbits much longer than $N=500$ steps. Our data sets are then built by restricting these long orbits to chunks $N$-steps long.
On the opposite, \cite{boulle2020classification} considers a ground truth labeling criterum based on a combination of Lyapunov exponent and Shannon entropy. 

In Appendix~\ref{app:opti}, we provide exhaustive analysis of the network performance considering various strides, activation functions and size of the first dense layers. As in~\cite{boulle2020classification}, we observe that increasing the number of convolutional layers does not improve the performance. This suggests that sequence information is maximally condensed by 2 convolutional layers.

\section{Uncovering the hidden workings of LKCNNs \label{sec:res}}

\subsection{LKCNNs can classify non-periodic orbits and outperforms Lyapunov reconstruction \label{sec:perform_LKCNN}}

We subdivide the input sequences corresponding to regular orbits in periodic and non-periodic. As in \cite{boulle2020classification} we train the LKCNN on the logistic log-spaced data set and test it on the sine-circle map, see Table \ref{tab:perform}. The training accuracy was set to 0.975. We observe low performance on the regular non-periodic subset. This is to be expected since the non-periodic subset of the logistic and sine-circle map are qualitatively different. If we train the network on the sine-circle data set we observe that the accuracy on the training set will not exceed 0.8. If we remove the periodic trajectories from the data set we observe accuracy exceeding 0.8, see Table \ref{tab:perform}. Hence, this is a labeling problem where the regular trajectories need to be subdivided into periodic and non-periodic trajectories to distinguish them from chaotic trajectories.


\begingroup
\begin{table}
\centering
\small
 \begin{tabular}{llccc}
\toprule
& & \multicolumn{3}{c}{Accuracy} \\
 \cmidrule(lr){3-5}
Training set & Test set & Chaotic $\mu \pm \sigma$ & \multicolumn{2}{c}{Regular $\mu \pm \sigma$} \\
 \cmidrule(lr){4-5}
& &  & Periodic  &  Non-periodic \\
\midrule
Logistic & Logistic  & 0.98 $\pm$ 0.0075 & 0.98 $\pm$ 0.0042 & 0.94 $\pm$ 0.014 \\
\quad log-spaced & \quad log-spaced & & & \\ 
Logistic & sine-circle  & 0.88 $\pm$ 0.090  & 0.93 $\pm$ 0.052 & 0.045 $\pm$ 0.01\\
\quad log-spaced & \quad linear-spaced & & &  \\
\bottomrule
\end{tabular}
\caption{Performance of the LKCNN on our data subsets. We observe that training on the logistic set generalizes well to sine-circle set on the chaotic and periodic data subset of the sine-circle map but performs poorly on the regular periodic set of the sine-circle map. The latter is to be expected since the non-periodic trajectories are qualitatively different. \label{tab:perform}}
\end{table}
\endgroup

\begingroup
\begin{table}
\centering
\small
 \begin{tabular}{llccc}
\toprule
Model &  $P_{\rm chaotic}$ & $P_{\rm regular}$ (non-periodic) \\
\midrule
LKCNN & 1.0 & 1.0 \\
Short-time Lyapunov exponent & 1.0  & 0.99 \\
\quad input sequence & & \\
Direct reconstruction Lyapunov exponent & 1.0 & 0.93 \\ 
\quad with Savitsky-Golay  & &  \\
\bottomrule
\end{tabular}
\caption{Precision of models trained on non-periodic sine-circle orbits (cf. Equation~\eqref{eq:precision}): We observe that the short-time Lyapunov exponent which is the Lyapunov exponent computed over the input sequence very accurately predicts the long time Lyapunov exponent. Hence, we can compare LKCNN prediction to prediction based on direct reconstruction of Lyapunov exponents via Savitsky-Golay polynomial filtering. The LKCNN outperforms the reconstruction.  \label{tab:diff_models}}
\end{table}
\endgroup

 In Table \ref{tab:diff_models}, we investigate the performance on non-periodic data for LKCNN, average Lyapunov exponent over the input sequence which we refer to as short time Lyapunov exponent and direct reconstruction of Lyapunov exponent with Savitsky-polynomial filtering, Section \ref{sec:SG}. 

Short time Lyapunov exponent is obtained with knowledge of $f$ and direct reconstruction with Savitsky-Golay is obtained without knowledge of $f$.  As the short time Lyapunov exponent is nearly perfect it implies that if the direct reconstruction can perfectly determine the derivatives then it can nearly perfectly classify the sequence. Hence, we shall compare the performance of the direct reconstruction to LKCNN. We observe that LKCNN performs significantly better.    
To add more rigor the data set contains sequences where the Lyapunov exponent has converged to $k$-decimal precision which are then rounded to $k$-decimals to determine the label. This methodology applies to the short time Lyapunov exponent and the direct reconstruction of Lyapunov exponent in Table \ref{tab:diff_models}. The results in Table \ref{tab:diff_models} are for $k=4$ but these results persist, see Appendix \ref{app:perform}.

\subsection{Mapping periodic input to periodic activation is model-independent over generic untrained LKCNN\label{sec:pres_per}}

\begin{figure}[p]
     \centering
     \begin{subfigure}{.75\textwidth}
         \centering
         \includegraphics[width=10cm]{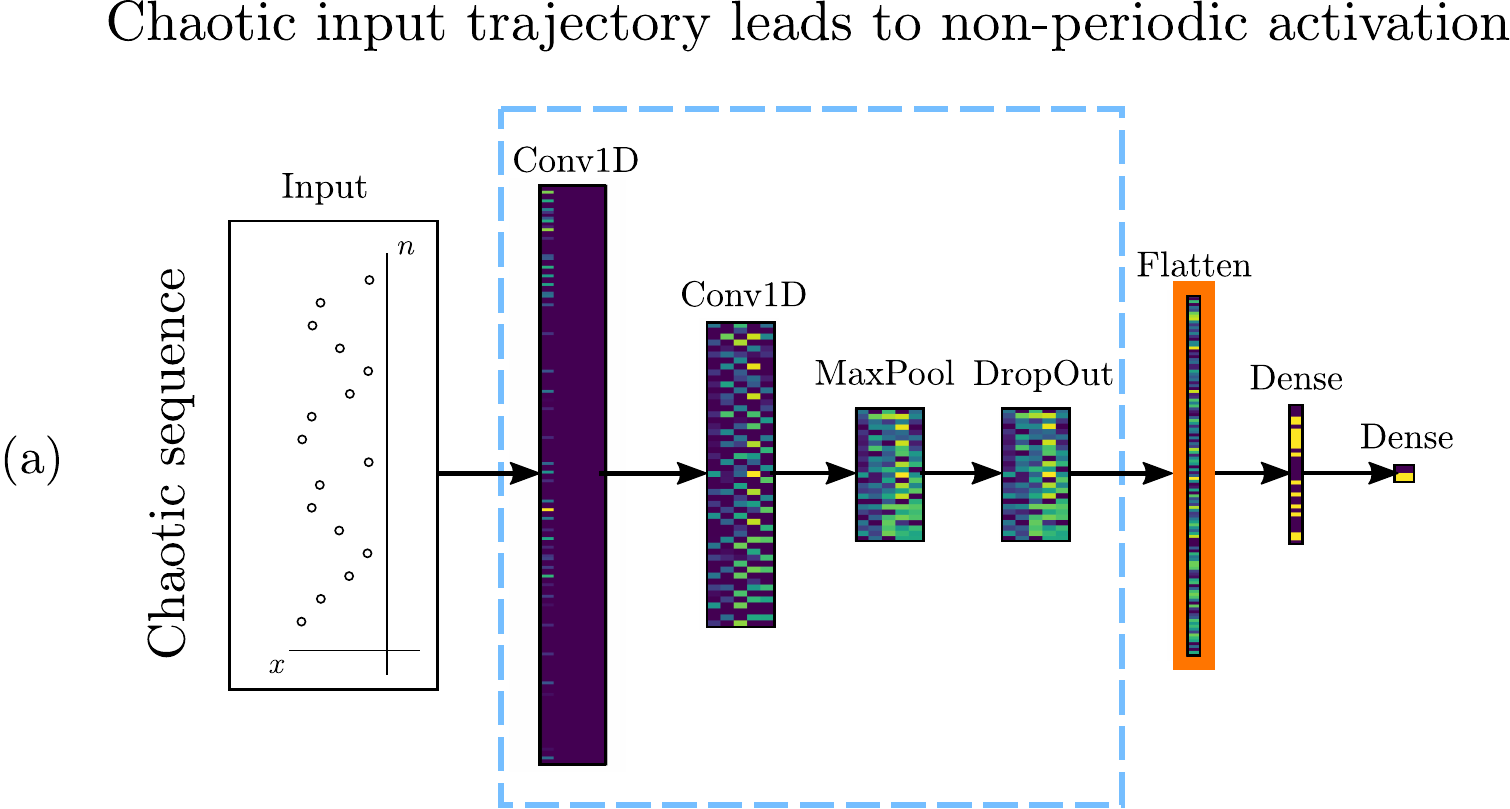}
         \caption*{ \label{fig:chaotic_activ}} 
     \end{subfigure}
    \begin{subfigure}{.75\textwidth}
         \centering
         \includegraphics[width=10cm]{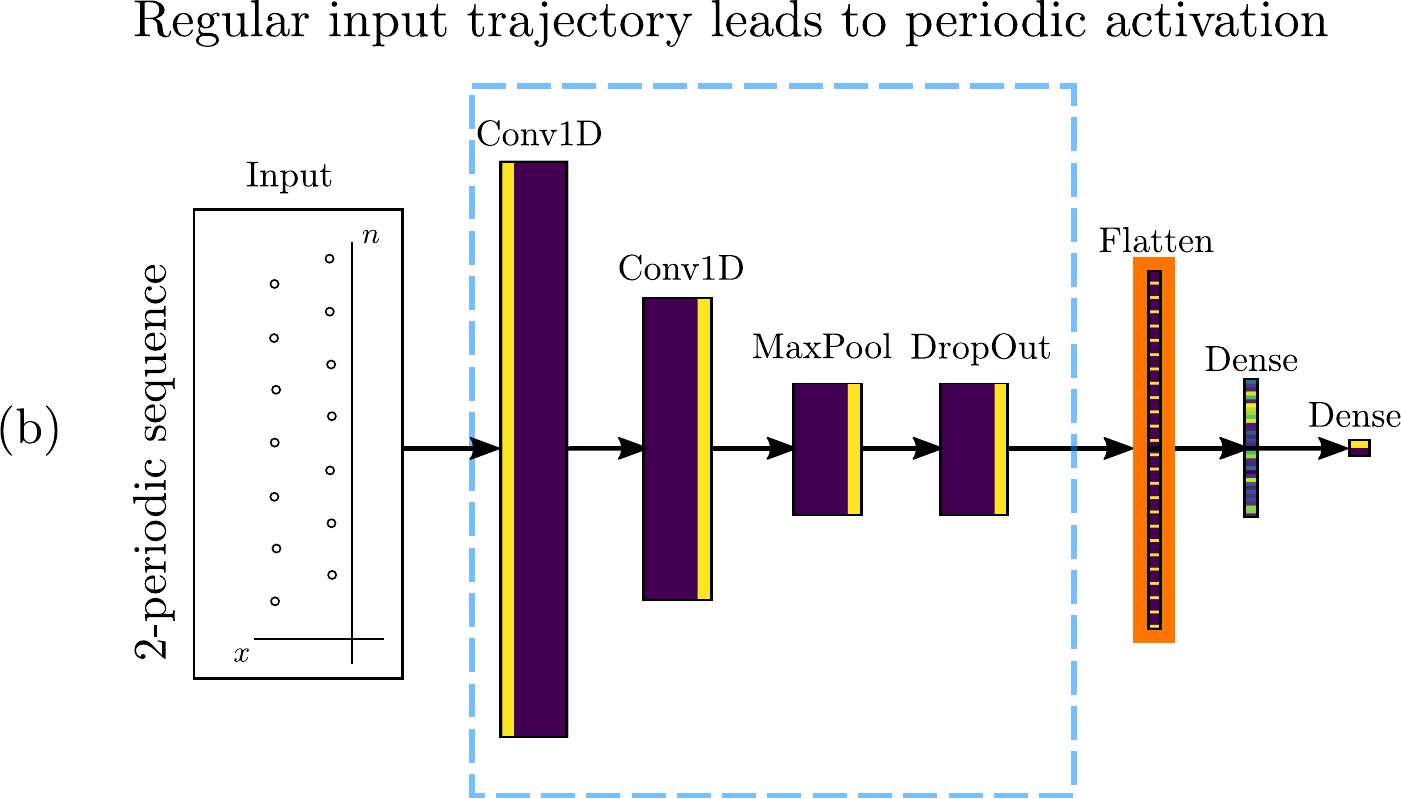}
         \caption*{  \label{fig:periodic_activ}}
     \end{subfigure}
    \caption{Activations LKCNN: We visualize the activation for a correctly classified chaotic and periodic trajectory. From the first convolutional layer to the flatten layer the non-periodic, (a), or periodic structure, (b), is preserved. The last two layers are dense and do not preserve any periodic or non-periodic structure.   \label{fig:all_activ}}    
\end{figure}

We consider the LKCNN trained on the log-spaced logistic data set, Section \ref{sec:data}. The activations for the first convolutional layer to the flatten layer preserves non-periodic and periodic structures. For a chaotic sequence the activations are non-periodic and for a periodic sequence the activations are periodic, see Figure \ref{fig:all_activ}.  We first present a rigorous result on the periodicity of a convolutational layer with stride $s$ for periodic inputs and then extend this result to LKCNN.  

\subsubsection{Periodic activations of convolutional layers \label{sec:rigor}}

Borrowing from equivariant theory for convolutional networks~\cite{cohen2016group}, we show that in an idealized setting where the network's input sequence has infinite length the convolutional filter applied to a $k$-periodic sequence yields an activation which is $k$-periodic. Furthermore, for stride $s$ the period becomes $k/s$ if $k$ is divisible by $s$. \\

Denote by $F: \Z \rightarrow \R$ the sequence feature map which is the map associated to a sequence $\{ x_{i}\}$ such that $F(i)=x_i$. Denote the convolutional filters by $\phi:\Z \rightarrow \R$. The convolutional filters also take as input the sequences indexes. A map $F$ is $k$-periodic if $F(z+k)=F(z)$ for all $z\in \Z$. Convolution of a sequence feature map $F$ by the filter $\phi$ is defined by
\begin{align*}
    [ F * \phi](z):= \sum_{y \in \Z} F(y) \phi(z-y)  .
\end{align*}
Define $S: \mathbb{Z} \rightarrow \mathbb{Z}$ by $S(z) = s z$ with $s \in \N$. Then a convolution with stride $s$ is given by  $[ F * \phi] \circ S$. Denote $k \in \N$ divisible by $s$ as $s|k$. 
\begin{lem}
Let $F$ be $k$-periodic. If $s | k$ then the period of $[ F * \phi] \circ S$ is $k/s$. If $s \nmid k$ then the period of $[ F * \phi] \circ S$ is $k$.  
\label{lem:conv_per}
\end{lem}
\noindent Observe that if $F$ is $k$-periodic and $s=1$ then the period of $[ F * \phi] \circ S$ is $k$. Note that $[ F * \phi]$ can have period less than $k$ for suitably chosen $\phi$. For example, if $\phi$ is the zero function then $[ F * \phi]$ has period 1. 

\begin{proof} Take $\hat{k} \in \N$. We define $\hat{y}= y-s\hat{k}$. We can write 
\[
F(y) \phi(s(z+\hat{k}) - y)  = F(\hat{y}-s\hat{k}) \phi(sz-\hat{y}).
\]
Using the above we can write
\begin{align*}
    [ F * \phi]S(z+\hat{k}) &= \sum_{ \hat{y} \in \Z } F(\hat{y}-s\hat{k}) \phi(sz-\hat{y}).
\end{align*}
Hence,  if $s | k$ then we obtain that $[ F * \phi]S(z+\hat{k}) = [ F * \phi]S(z)$ for  $\hat{k} = k/s$ and if 
$s \nmid k$ then we obtain that $[ F * \phi]S(z+\hat{k}) = [ F * \phi]S(z)$ for  $\hat{k} = k$.
\end{proof}

\subsubsection{Periodic activations of LKCNN \label{sec:per_LKCNN}}

The result from Section \ref{sec:rigor} assumes an infinitely large network. Practically, the periodicities that can be captured by the network depend on the size of the network.  Observe that the largest activation period that the first convolutional layer can capture is 100 since it has size $201 \times 5$, see Figure \ref{fig:all_activ}. Hence, since we have stride 2 the largest period the network can capture  in terms of the input $x$ is 200. 

Denote by $k$ the period of the activation for Maxpooling and $p$ the pool size. Then if $p | k$ the output period is $k/p$ and if $p \nmid k$ then the output period is $k$. Here, we have $p=2$. Again, in a practical setting the size of the layer restricts the periods that can be captured. Here, the output activation period can be at most 12. 

\begin{figure}[ht]
\centering
\includegraphics[width=10cm]{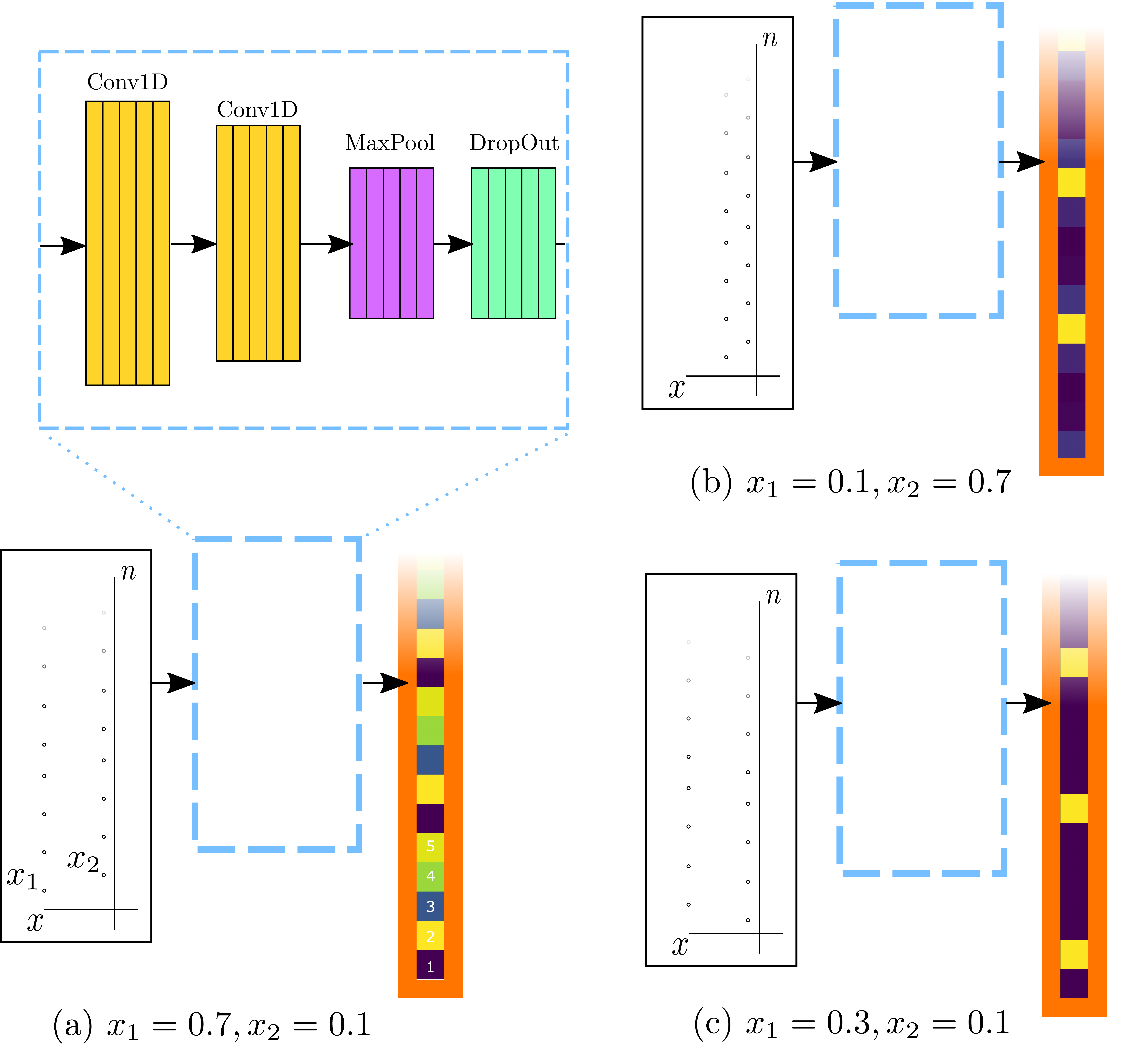}
\caption{Activation at flatten layer: We consider the first 15 nodes at the flatten layer for period 2 trajectories given by $x_1,x_2$. We observe that the periodicity is invariant for all three cases. Reversing the order alters the activation (a)(b). Position of the maximum activation, yellow, and minimal activation, purple, can vary (a)(b)(c). \label{fig:period2}}  
\end{figure}

The flatten layer will have periodic activation if the dropout layer also has periodic activation. Furthermore, if the activation matrix at the dropout layer is non-constant then the periodicity is increased by a multiplicative factor of 5 since the number of columns of the activation matrix is 5 which is prime.  Our experiments indicate that activation at the flatten layer can vary if we vary over the period of $x$, Figure \ref{fig:period2} and Appendix \ref{app:gen-2-per}.

The activation periodicity is lost in the last two dense layers. Hence, to study network periodicity we will focus on the periodicity at the flatten layer in relation to the periodicity of the sequence.

Combining the results from this section with Lemma \ref{lem:conv_per} we obtain the following heuristic statement: \textit{for an untrained LKCNN we expect that for $k\leq 96$ the period at the flatten layer is $ 5 k/2^i$ where $i$ is the largest $i \leq 3$ for which  $2^i|k$}. 


\subsection{Activation periodicity influences performance}

\begin{figure}[ht]
    \centering
 \includegraphics[width=12cm]{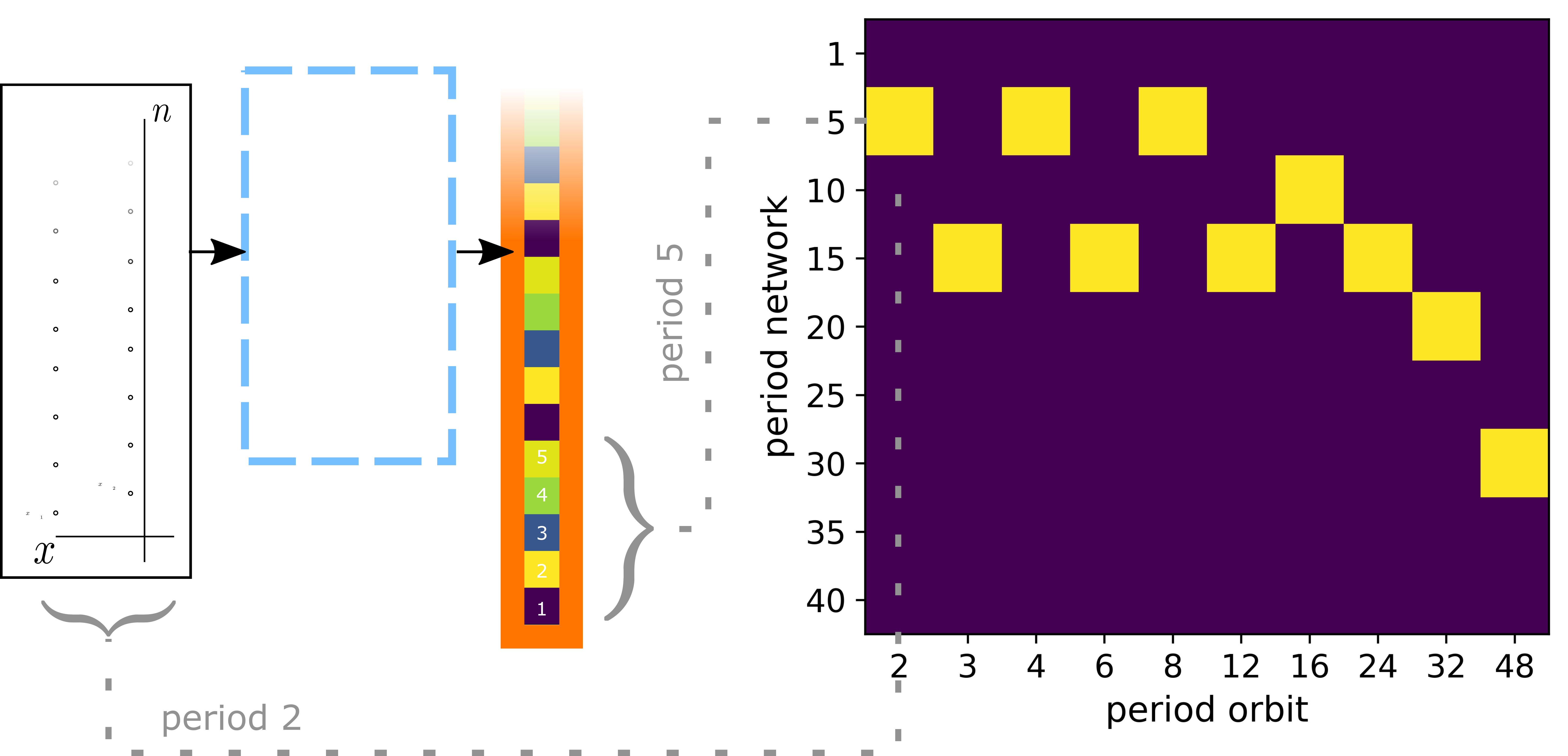}
\caption{Period matrix: On the left a period 2 orbit with flatten layer period 5 is depicted. For each periodic orbit we identify the orbit's period to a period at the flatten layer. The period at the flatten layer is referred to as network period. Orbits with the same period do not uniquely map to a network period. However, in period range 2 to 48 the majority of period orbits map to a single network period. We associate a binary matrix to this identification. This period matrix can differ per trained model. However, we can identify classes of models which have the same period matrix \label{fig:periodA}}
\end{figure}

Generally, a periodic input implies a periodic activation at the flatten layer. We refer to the latter as the network period. For example, the 2 period orbit in Figure \ref{fig:period2} has network period 5.  Orbits with the same period map into a single network period if the orbit period is sufficiently small. Consequently, we represent this mapping by a binary matrix which identifies orbit periods to network periods, see Figure \ref{fig:periodA}. 
\begin{figure}[ht]
     \centering
     \begin{subfigure}{.45\textwidth}
         \centering
         \includegraphics[width=6cm]{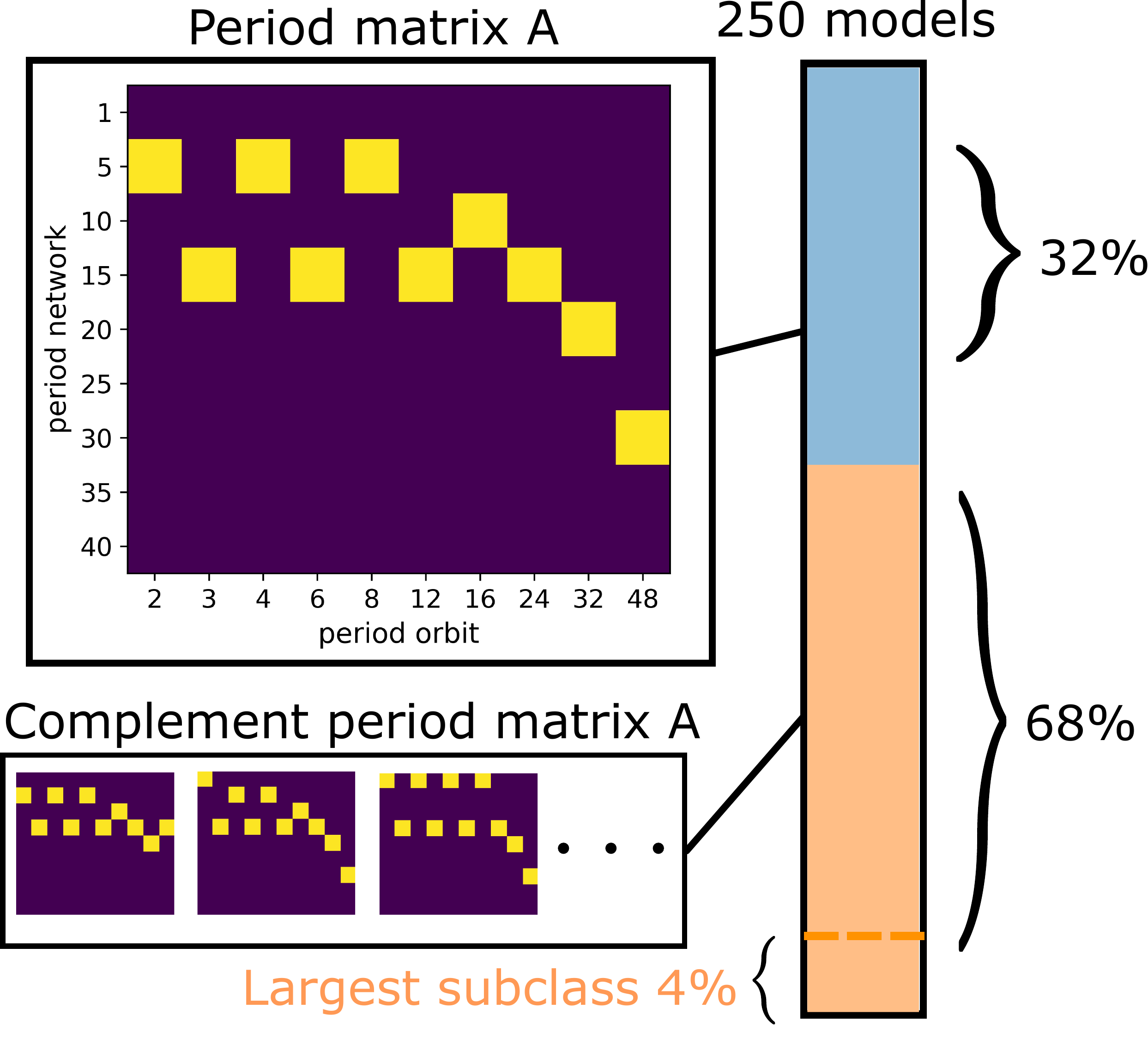}
         \caption{\label{fig:Aoverview}}
     \end{subfigure}
     \begin{subfigure}{.45\textwidth}
         \centering
         \includegraphics[width=4.5cm]{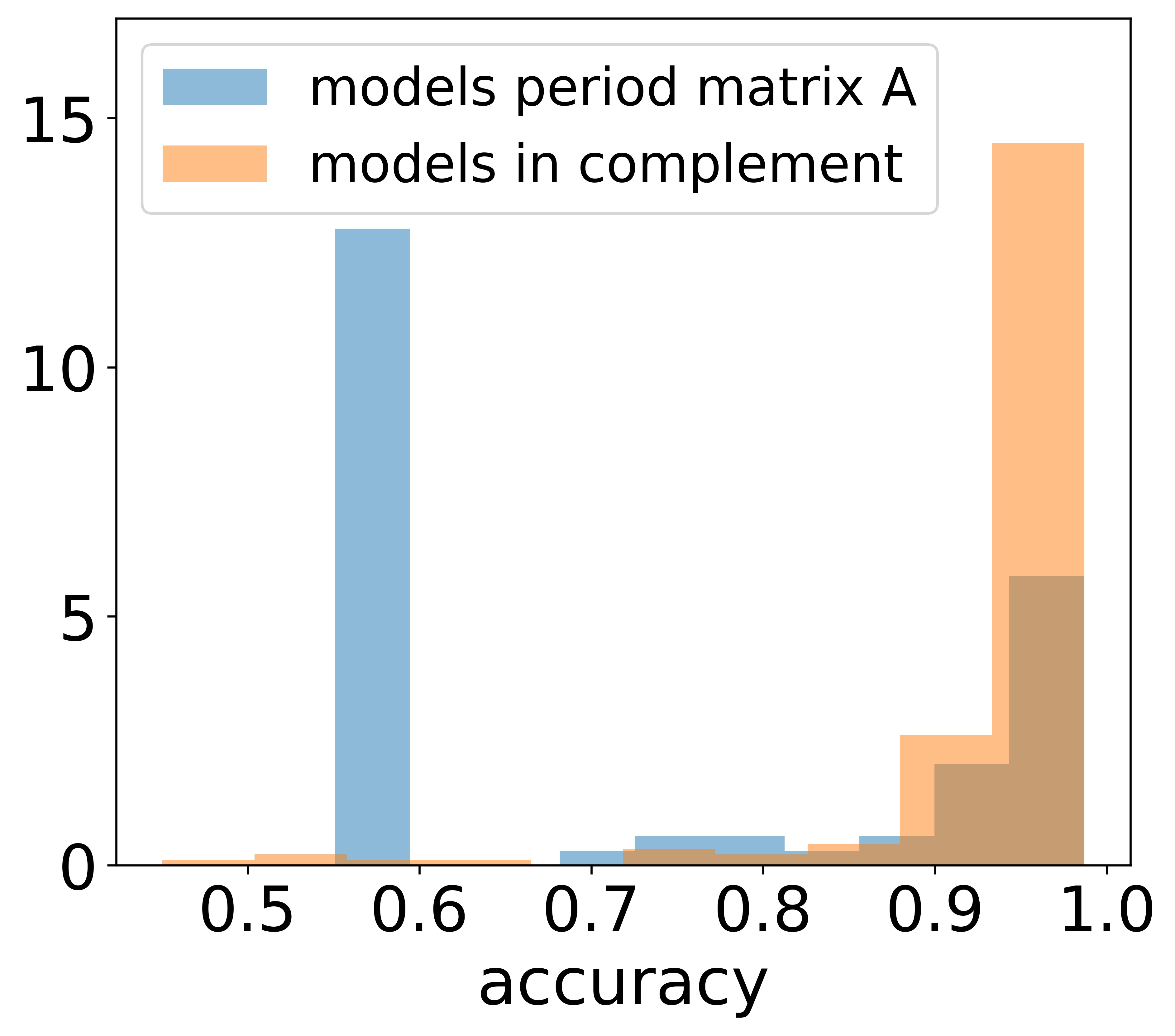}
         \caption{ \label{fig:norm_acc_hist}}
     \end{subfigure}
\caption{Models classified by period matrix: We investigated the period matrices associated to 250 random initialized trained models. In (a) we observe that the most common period matrix is period matrix $A$ to which 32\% of models converge. The second largest period matrix class is 4\%. We compare the performance between models with period matrix $A$ and models in their complement. In (b) we plot the density functions for the accuracy of both classes. We observe that models with period matrix $A$ have on average lower performance.  }     
\end{figure}
We note that we could formulate these results in a non-binary setting where the orbit period is not uniquely mapped to a single network period. However, similar results still hold, see Appendix \ref{app:non-bin}.

The binary matrix in Figure \ref{fig:periodA} is called the period matrix.  Period matrices can differ per trained model. We consider 250 models trained over 1000 epochs with patience 50. We then classify the trained models using the period matrices. The largest class consist of 32\% of trained models and has period matrix as depicted in \ref{fig:Aoverview} which will be referred to as period matrix $A$. Single classes in the complement of $A$ make up at most 4\% of all models. We consider these classes to be non-significant. We note that the great variety in period matrices is not in contradiction with Lemma \ref{lem:conv_per} since the Lemma only gives an upper bound on the activation period of the convolutional layers.

We observe that 56\% of class $A$ models converges to a local minimum which has accuracy between 55-60\%, Figure \ref{fig:norm_acc_hist}.  For the models in the complement of class $A$ we have that Q1 and Q3 are at 0.94 and 0.98 accuracy, respectively. Hence, models of class $A$ have a negative effect on the overall performance. 

The period matrix $A$ exactly corresponds to the theoretically determined periods for untrained networks, Section \ref{sec:per_LKCNN}. For other period matrices the weights of the corresponding models reduce the network period. This appears to generally have a positive effect on the performance. Heuristically, we can argue that there is no benefit to preserve a property of the untrained network which is imposed on the network independent of data properties. Therefore, models with period matrix $A$ underperform in relation to its complement. 

\section{Discussion \label{sec:conc}}
In this work we consider the issue of classifying  time sequences, discriminating regular from chaotic dynamics. To this purpose we compare and analyze a state-of-the-art approach based on Large Kernel Convolutional Neural Networks (LKCNNs) and a traditional Lyapunov reconstruction method.
LKCNNs have a simple structure: few convolutional layers with peculiarly large filters connected to a dense layer.  This structure simplifies numerical experiments and allowed us to identify relevant network features key for performance.
Specifically, we have shown that to classify signals with high accuracy,  LKCNNs use qualitative properties of the input sequence. This enables them to outperform direct  Lyapunov exponent  reconstruction methods. Here, we consider a reconstruction approach based on Savitsky-Golay polynomial filtering. LKCNNs higher classification accuracy strongly emerges as we consider sequences whose Lyapunov exponent is close to zero: the hardest to evaluate as the average rate of separation between consecutive time-steps is almost constant.  We investigated the emerging connection between input periodicity and periodicity in the activation of the non-dense layers. We have shown this aspect to be paramount for performance. For this analysis, we have introduced the notion of a period matrix which is a two dimensional binary matrix which represents the mapping between periodic inputs and the network's activation excluding the dense layers. 

We considered generic untrained networks first. For these, we have determined theoretically the connection between the input period and the periodicity of the activation of the convolutional layers. Effectively, this yields a period mapping-type relation. This mapping is represented by the period matrix. We showed that for generic networks this period matrix is unique. The weights minimizing the loss function are meager within the weight space. Consequently, this period mapping needs not be preserved by training. Indeed, we observed that trained models can have a variety of period mappings. Nevertheless, a significant percentage of trained models have period matrix equal to that of generic untrained models. These models underperform. In other terms, we numerically verified how the period matrix is a feature correlating with performance.
Heuristically, there is no benefit to preserve a property of the untrained network if it has no relation to the data. If a property of the data is reflected on the network we would expect an increase in performance. 

Models with high performance have period matrices featuring network periods which are lower compared to the case of a generic untrained network. Additionally, high performance is not a property of a singular period matrix. Two aspects remain outstanding: how the network training yields period reduction, and whether properties of the periodicity matrix correlating with high classification performance can be identified independently on the periodicity matrix of the untrained network. Possibly, the absence of clear connections indicates that restrictions on the period matrix adversely affect performance.

Finally, here we focused on discrete dynamical systems. It would be relevant to formulate suitable experiments in a continuous setting. We note that our approach leverages on periodicity of trajectories. Hence, it could be applied to other discrete dynamical systems or settings where it is possible to construct Poincar\'{e} maps. \\

This paper illustrates how dynamical system problems can be solved using neural networks and how neural network problems can be understood using dynamical systems. Thanks to a neural network approach we could achieve higher performance than traditional methods at classifying chaotic from regular time series. On the other hand, an approach hinged on dynamical systems analysis has been key to understand the inner workings of the network. \\


\noindent \textbf{Acknowledgments:} During this research Thomas de Jong was also affiliated to University of Groningen and Xiamen University. Many thanks to Alef Sterk for his helpful comments and literature recommendations. Also, many thanks to Klaas Huizenga for providing hardware during Thomas de Jong's stay in Groningen. This research was partially supported by JST CREST grant number JPMJCR2014.
\appendix

\part*{Appendix}

\section{Generalizations of 2-periods \label{app:gen-2-per}}

\begin{figure}[ht]
     \centering
     \begin{subfigure}{.45\textwidth}
         \centering
         \includegraphics[width=\linewidth]{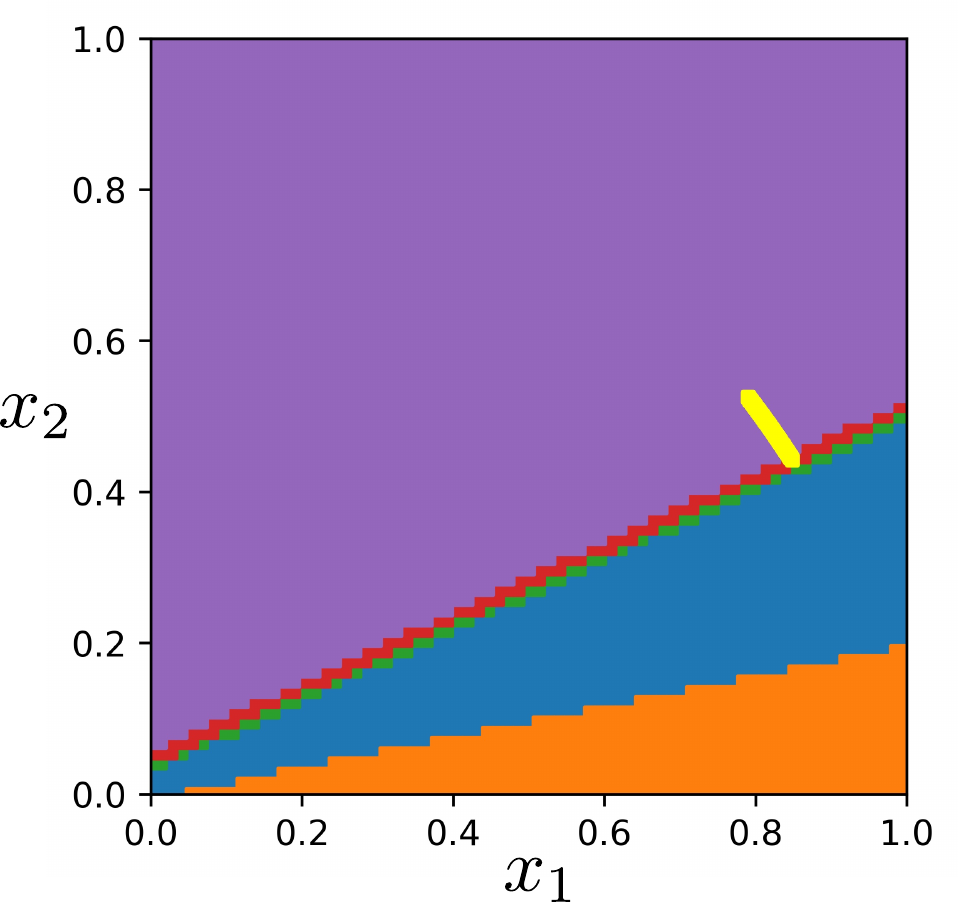}
         \caption{Bifurcation diagram model 1 \label{fig:bifu_model1}}
     \end{subfigure}
    \begin{subfigure}{.45\textwidth}
         \centering
         \includegraphics[width=\linewidth]{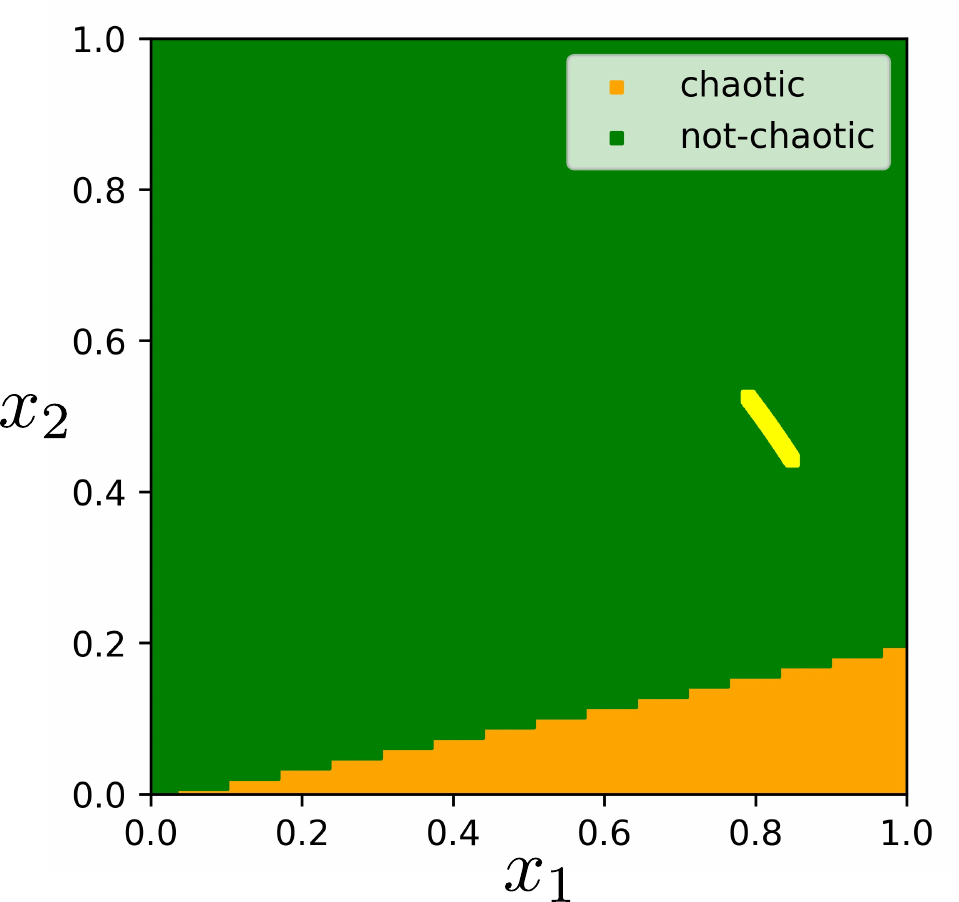}
         \caption{Prediction diagram model 1 \label{fig:pred_model1}}
     \end{subfigure}
     \begin{subfigure}{.45\textwidth}
         \centering
         \includegraphics[width=\linewidth]{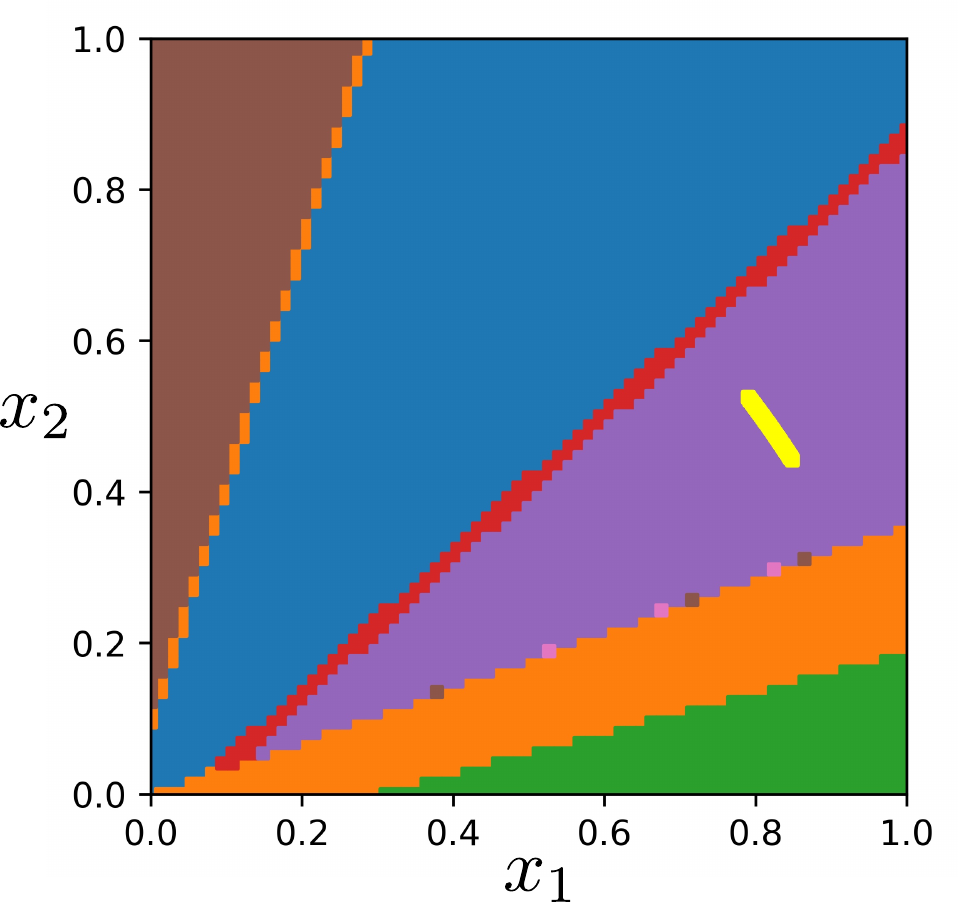}
         \caption{Bifurcation diagram model 2 \label{fig:bifu_model2}}
     \end{subfigure}
     \begin{subfigure}{.45\textwidth}
         \centering
         \includegraphics[width=\linewidth]{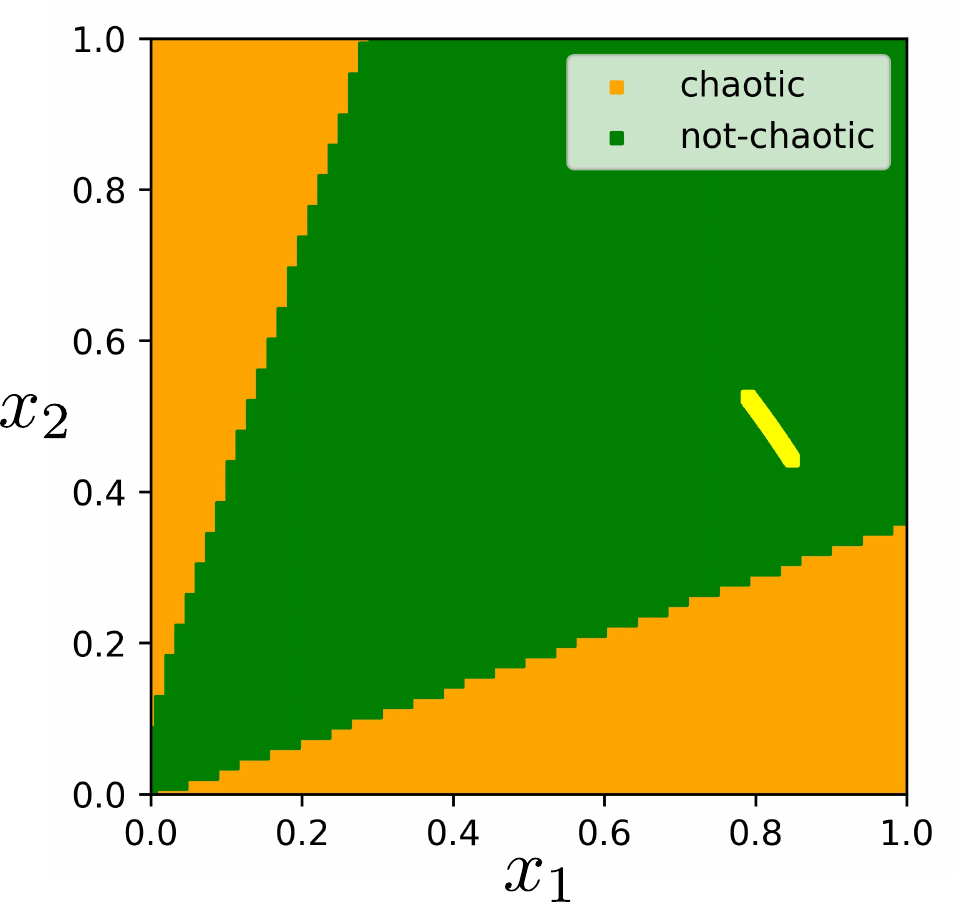}
         \caption{Prediction diagram model 2 \label{fig:pred_model2}}
     \end{subfigure}
\caption{Period-2 orbits bifurcation diagrams with predictions. The yellow line corresponds to period-2 orbits in the data. The other colors in (a)(b) correspond to period-2 orbits with zero-activation at the same position. \label{fig:bifu_diag}}     
\end{figure}

Using the method from Section \ref{sec:pres_per} we can identify the activation at the flatten layer for a single orbit. We can consider all possible orbits and investigate to which activation they are mapped. 
For visualization purposes we will only consider period-2 orbits. Recall from Figure \ref{fig:period_distribution} that the log-spaced data set does not contain fixed points. We consider two trained models. For both models the period at the flatten layer is either period 1 or 5.  This results in the bifurcation diagrams \ref{fig:bifu_diag}. 

In Figure \ref{fig:bifu_diag} the first and second iterate of the period-2 orbit are on the $x$-axis and $y$-axis, respectively. We have identified period-2 orbits by the position of their zero activations at the flatten layer, e.g. an activation of $[0,0,0,1,0]$ has the same color as an activation of $[0,0,0,0.1,0]$ as their zeros have the same position.

\begin{figure}[h]
    \centering
    \includegraphics[width=6cm]{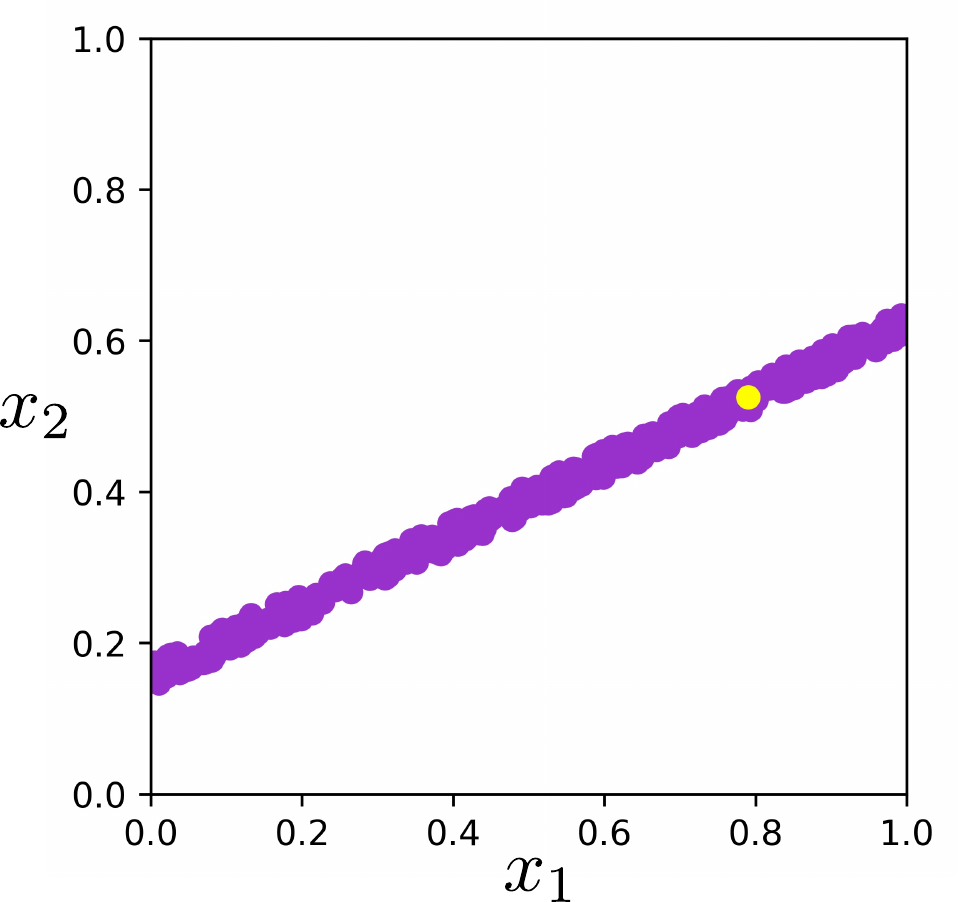}
    \caption{Sub-domain of Figure \ref{fig:bifu_model1}: The purple domain corresponds to a sub-domain of the purple domain in Figure \ref{fig:bifu_model1}. This domain has activation close to the yellow data point in the 2-norm.}
    \label{fig:p2act}
\end{figure}

We observe that for all diagrams in Figure \ref{fig:bifu_diag} the domains are bounded by linear functions of the form $x_2= a x_1 + b$ with $a >0$. Specifically, if we consider the activations which are close to a period-2 orbit in the 2-norm we also obtain these linear functions, see Figure \ref{fig:p2act}.  

In the bifurcations diagrams in Figure \ref{fig:bifu_model1} and Figure \ref{fig:bifu_model2} the activation classification is not preserved under changing the order of $x_1,x_2$. This order preservation is also not satisfied for the prediction, see Figure \ref{fig:pred_model1} and Figure \ref{fig:pred_model2}. As expected the prediction is incorrect if the prediction is applied for period-2 orbits far away from the trained data. In Figure \ref{fig:bifu_model2} there are 9 domains (excluding the yellow data domain). This means that there are different permutations of zero activation with the same number of zeros since the period at the flatten layer is at most 5. From Figure \ref{fig:bifu_model1} we observe that this activation classification is not locally preserved in a neighborhood of the data since the yellow line borders the purple and red domain.

\qquad \newpage

\section{Architecture optimization \label{app:opti}}

Our experiments show that the performance can be improved by adjusting the following two hyperparameters:
\begin{itemize}
    \item[-] Stride in the two convolutional layers,
    \item[-] Number of nodes in the fully connected dense layer after the flatten layer.
\end{itemize}

All hyperparameter variations exhibited instability when trained using the log-spaced data set. We resolved this by setting the learning rate to $0.000388$.

\begin{figure}[h]
    \centering
    \includegraphics[width=12cm]{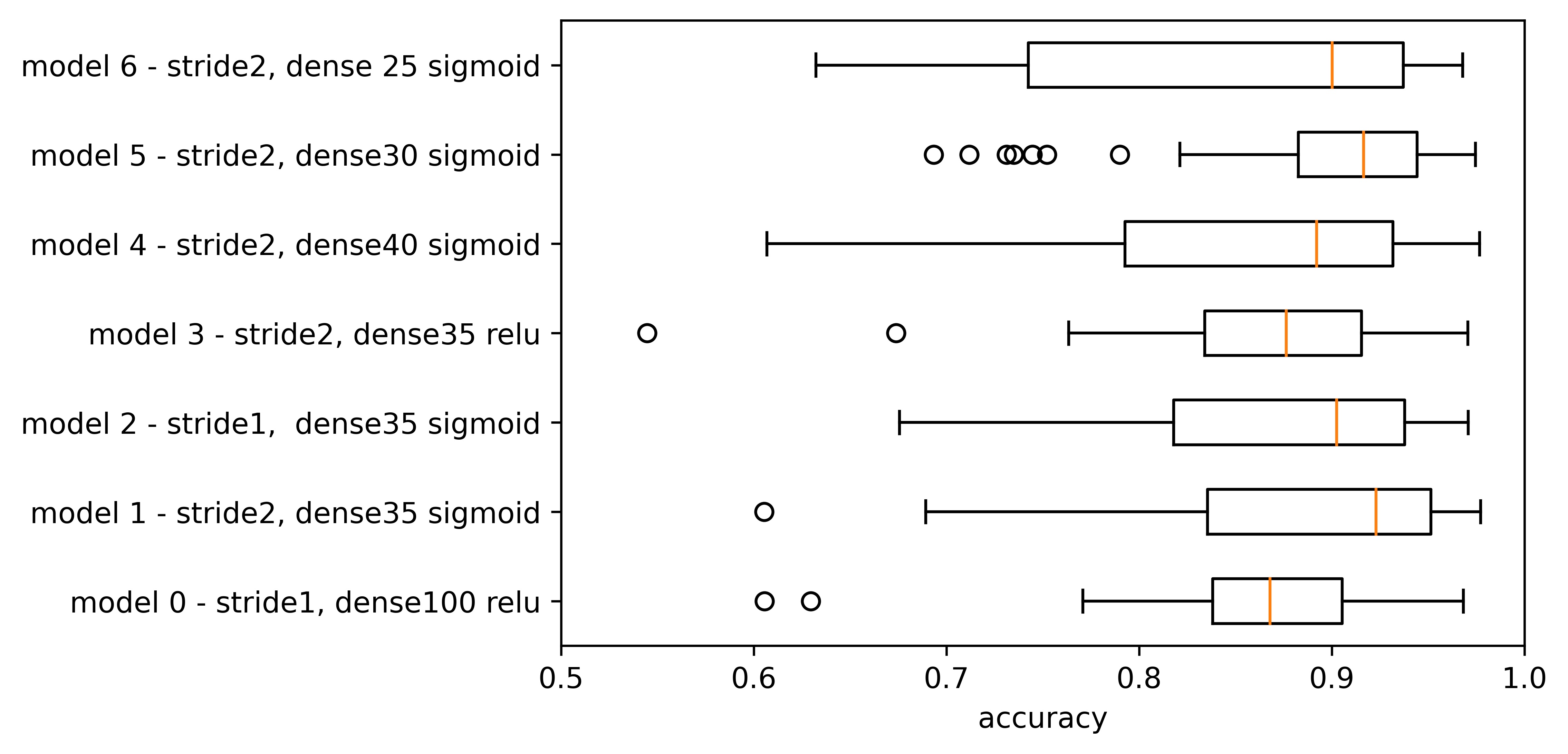}
    \caption{Validation accuracy boxplots for different hyperparameters. }
     \label{fig:acc_models}
\end{figure}

The experiments in Figure \ref{fig:acc_models} concern the accuracy over the validation set of 50 models for 200 epochs per architecture. Model 0 is the LCKNN architecture used in~\cite{boulle2020classification}. Since most of the experiments concerns models trained with high accuracy we are interested in selecting the hyperparameters which can consistently train models with high accuracy. We selected model 1 since the accuracy range over Q2 to Q3 is the highest.

\newpage
\section{Period matrices \label{app:non-bin}}

\subsection{Period matrix classes}
\begin{figure}[h]
     \centering
    \begin{subfigure}{.45\textwidth}
         \centering
         \includegraphics[width=\linewidth]{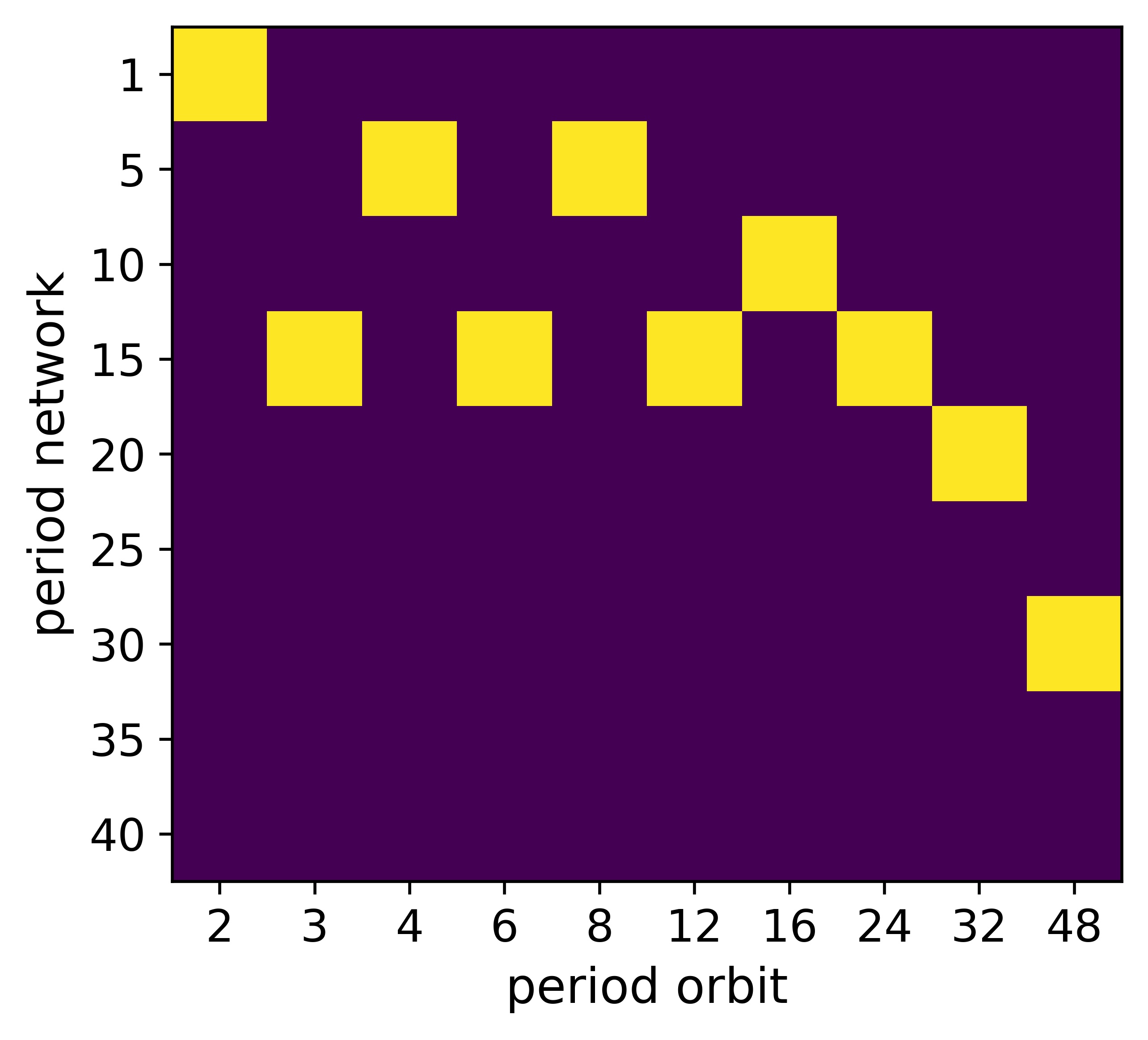}
         \caption{4\% \label{fig:11max}}
     \end{subfigure}
     \begin{subfigure}{.45\textwidth}
         \centering
         \includegraphics[width=\linewidth]{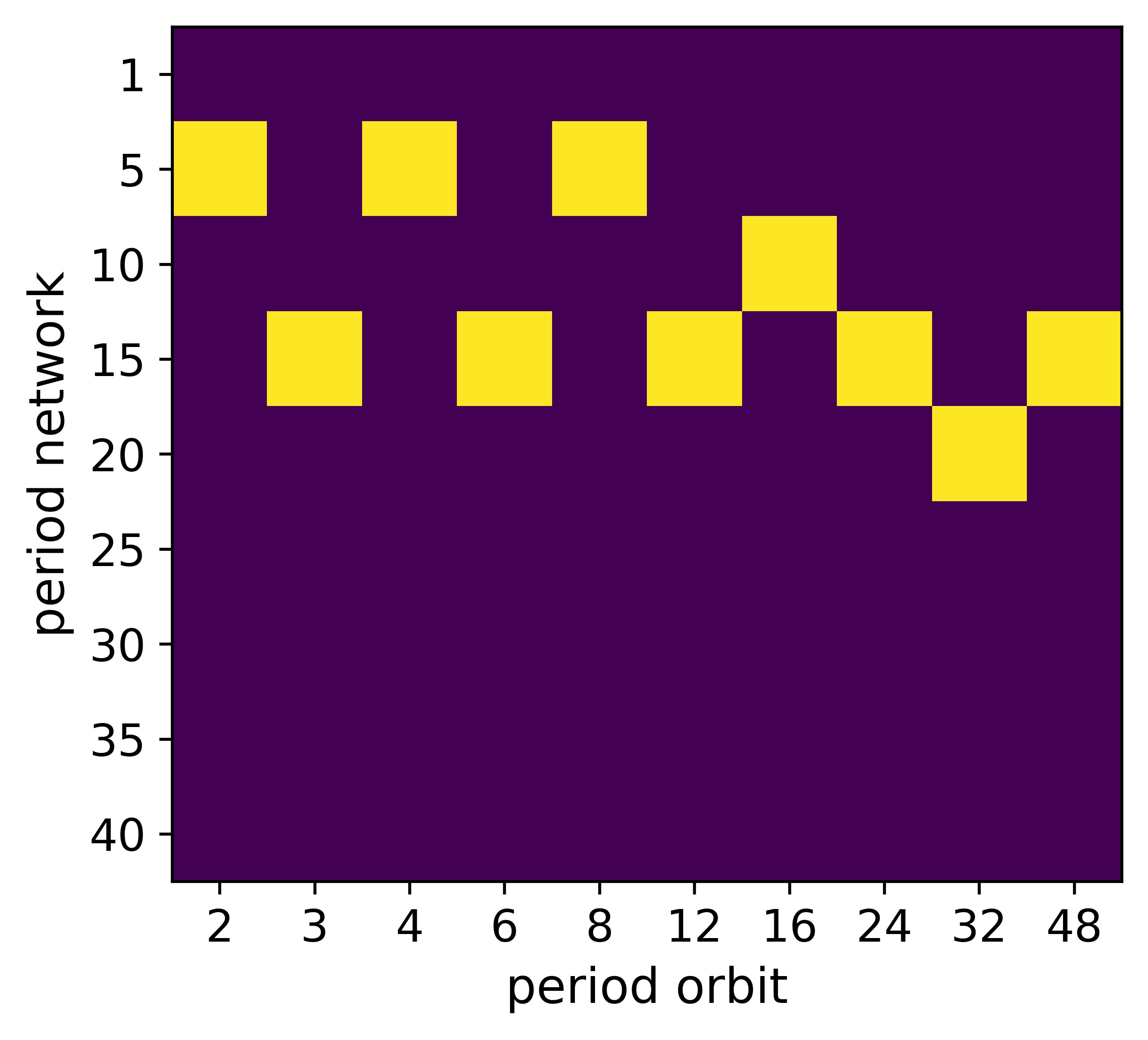}
         \caption{ 4\% \label{fig:7max}}
     \end{subfigure}
     \begin{subfigure}{.45\textwidth}
         \centering
         \includegraphics[width=\linewidth]{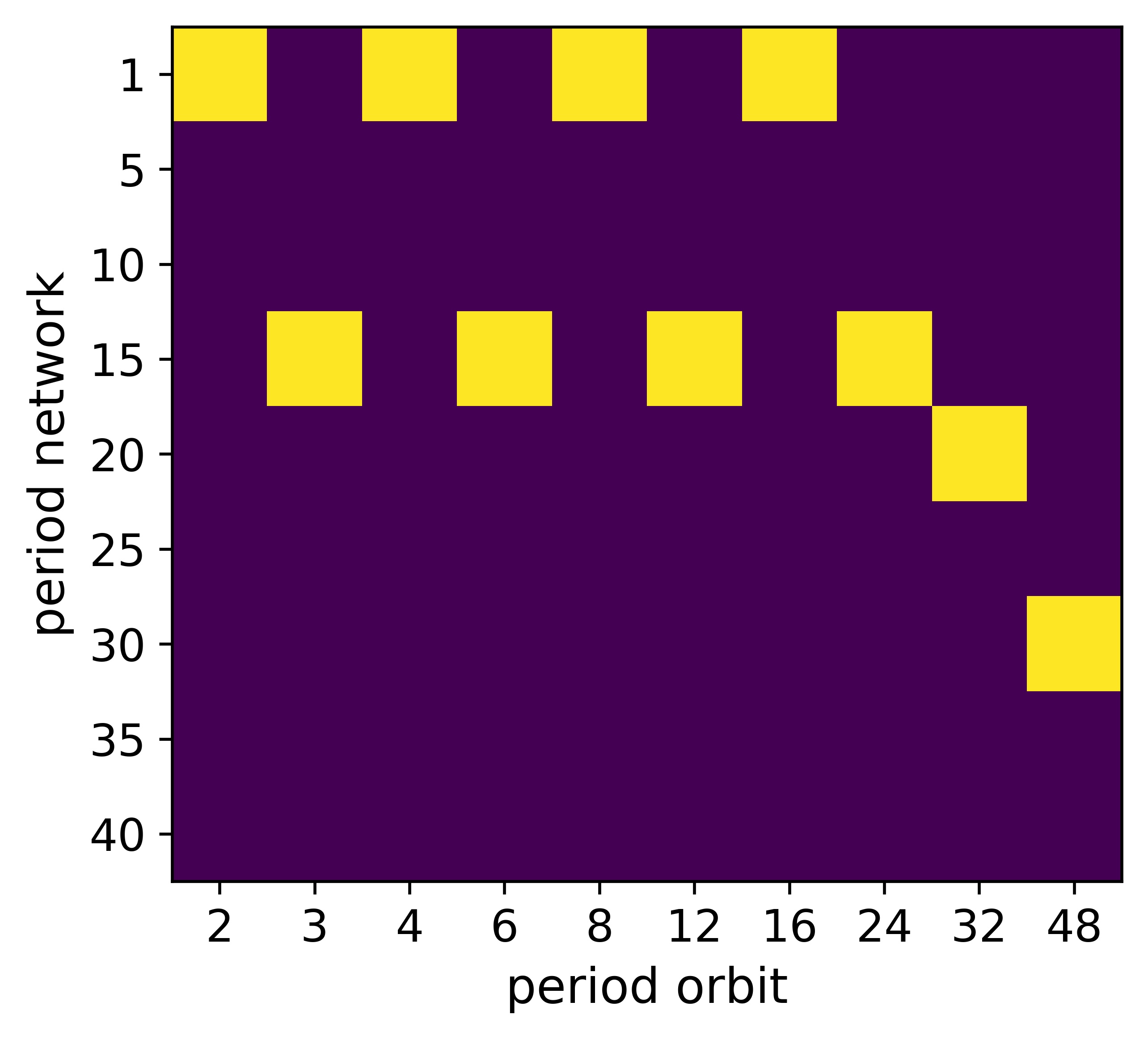}
         \caption{ 3\%  \label{fig:14max}}
     \end{subfigure}
\caption{Period matrices: The percentage indicates the size of the class with respect to all the models\label{fig:maxperiod}}     
\end{figure}

We recall that period matrix class A, Figure \ref{fig:Aoverview} makes up 32\% of the 250 trained models. In Figure \ref{fig:maxperiod} we visualize the 2nd, 3rd and 4th largest classes.

\subsection{Non-binary period matrices}

\begin{figure}[h]
     \centering
     \begin{subfigure}{.45\textwidth}
         \centering
         \includegraphics[width=\linewidth]{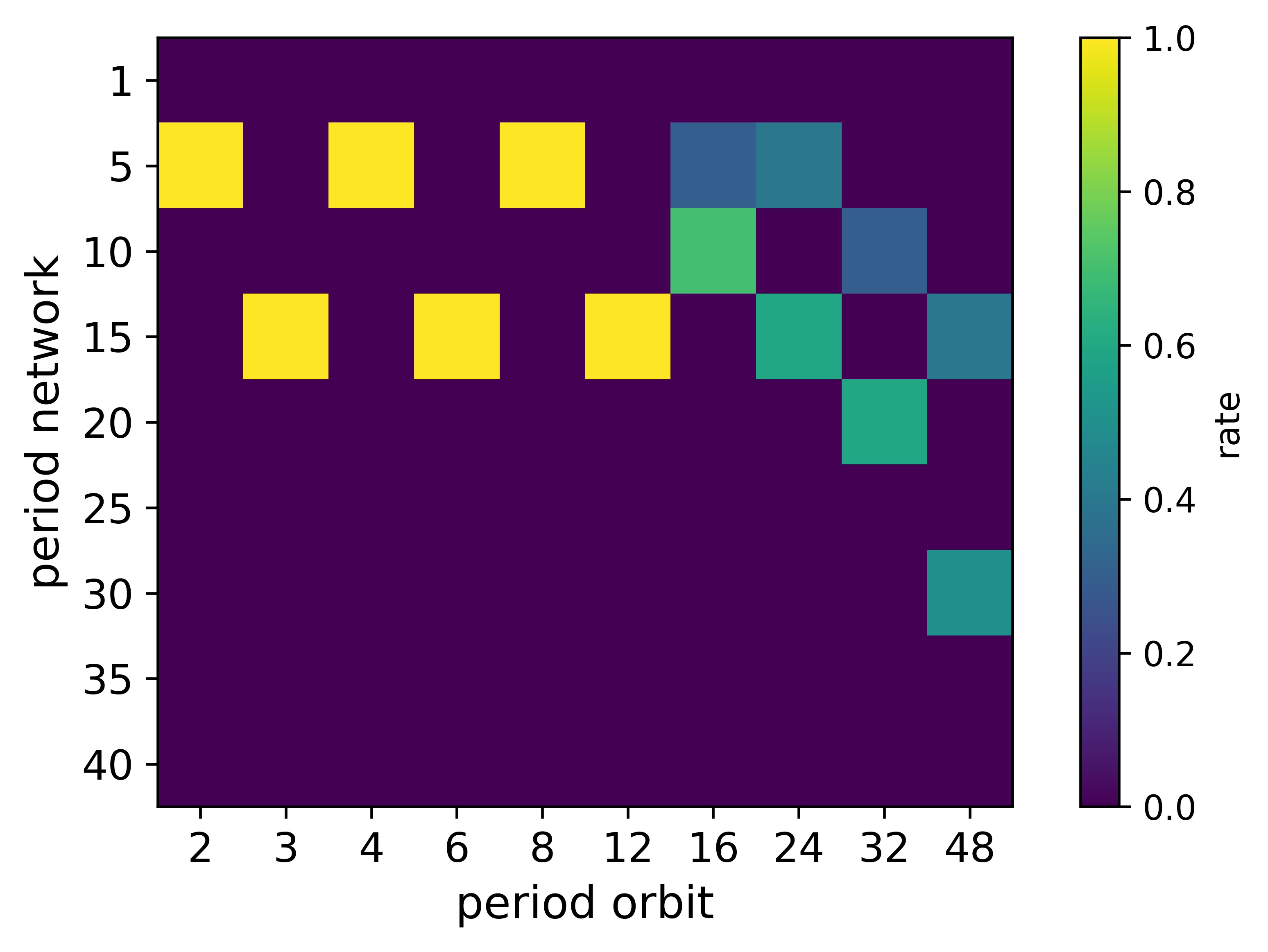}
         \caption{  29\% \label{fig:rawA}}
     \end{subfigure}
    \begin{subfigure}{.45\textwidth}
         \centering
         \includegraphics[width=\linewidth]{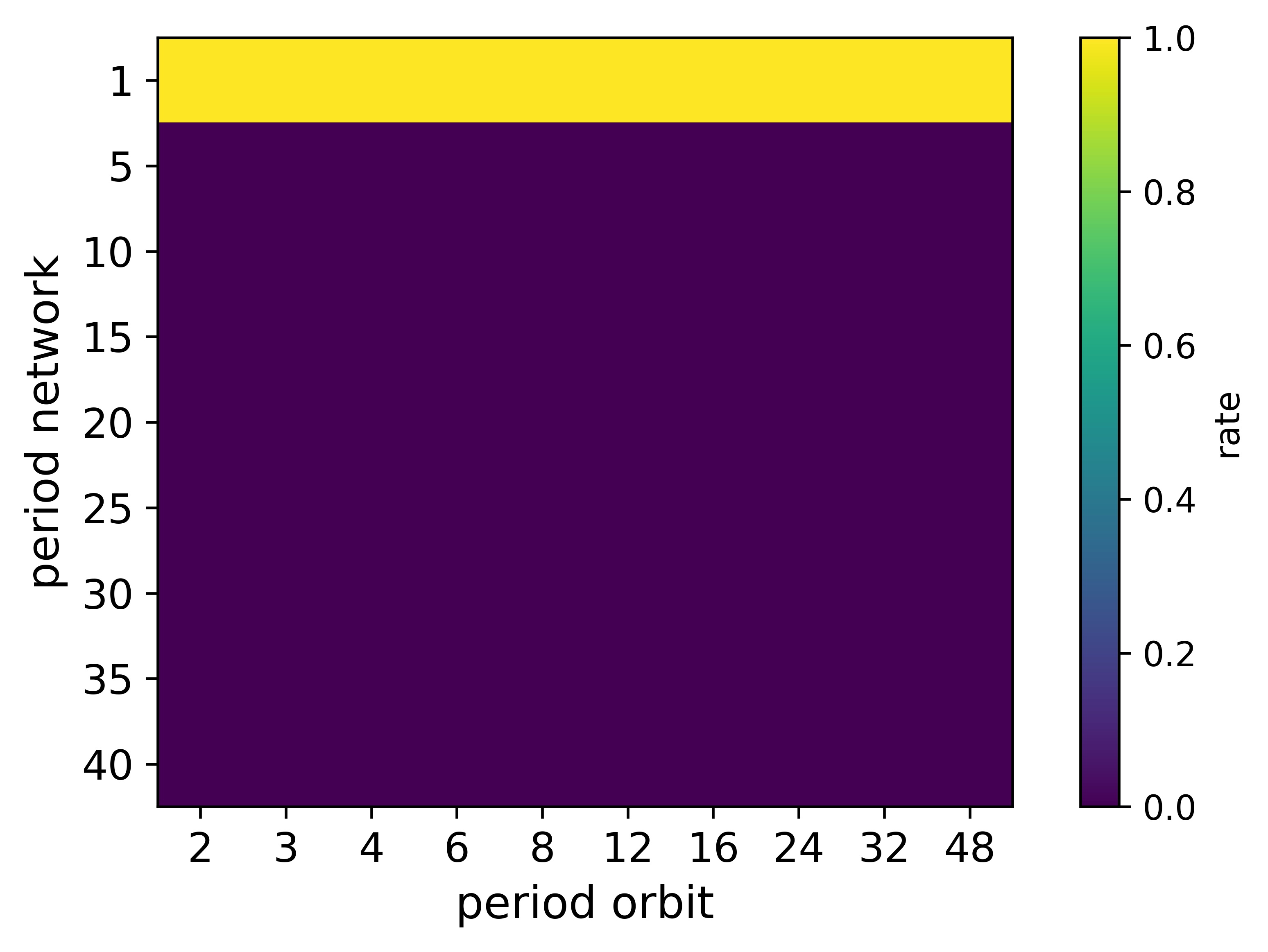}
         \caption{ 2\% \label{fig:52raw}}
     \end{subfigure}
     \begin{subfigure}{.45\textwidth}
         \centering
         \includegraphics[width=\linewidth]{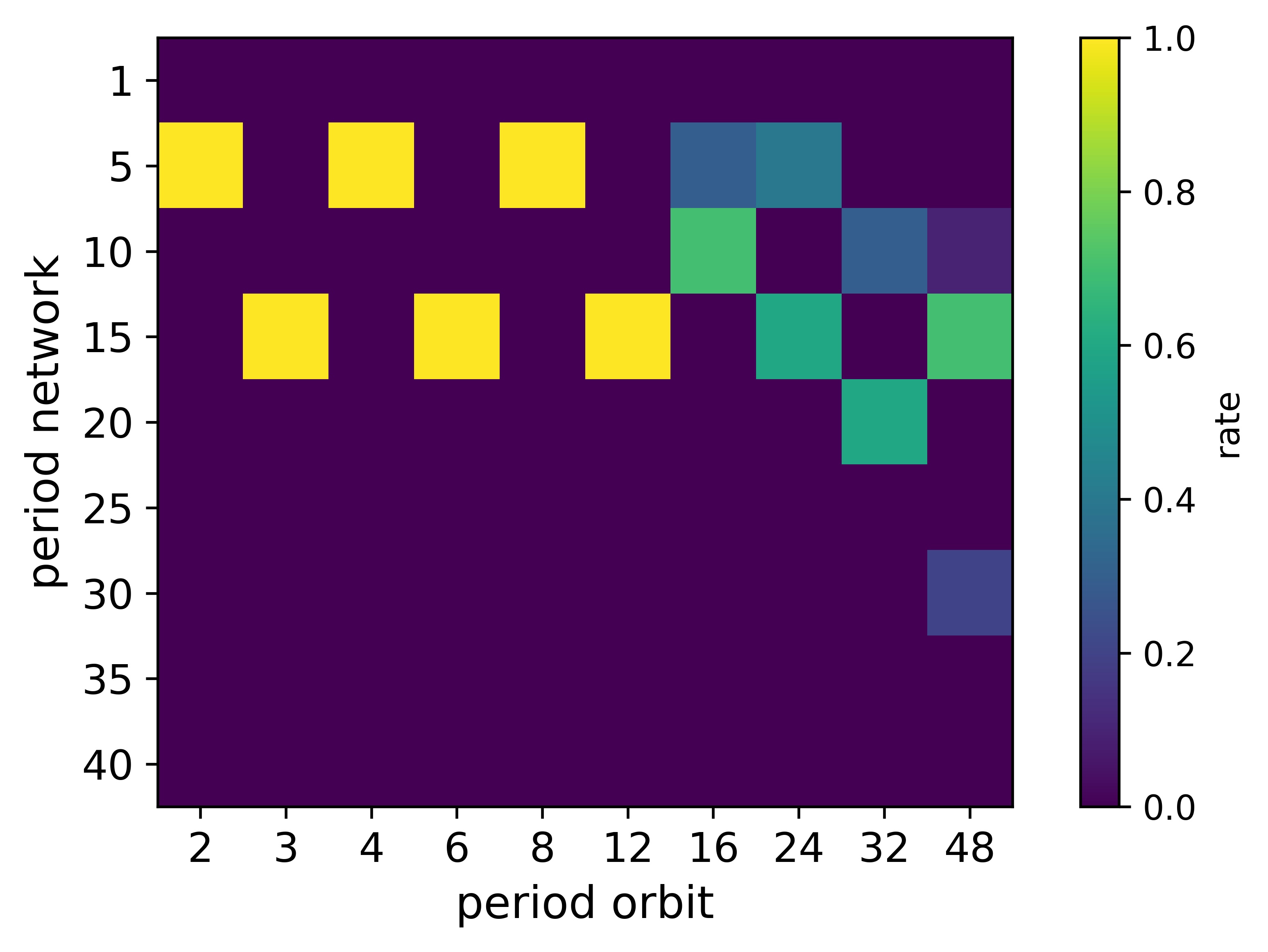}
         \caption{ 1\% \label{fig:49raw}}
     \end{subfigure}
     \begin{subfigure}{.45\textwidth}
         \centering
         \includegraphics[width=\linewidth]{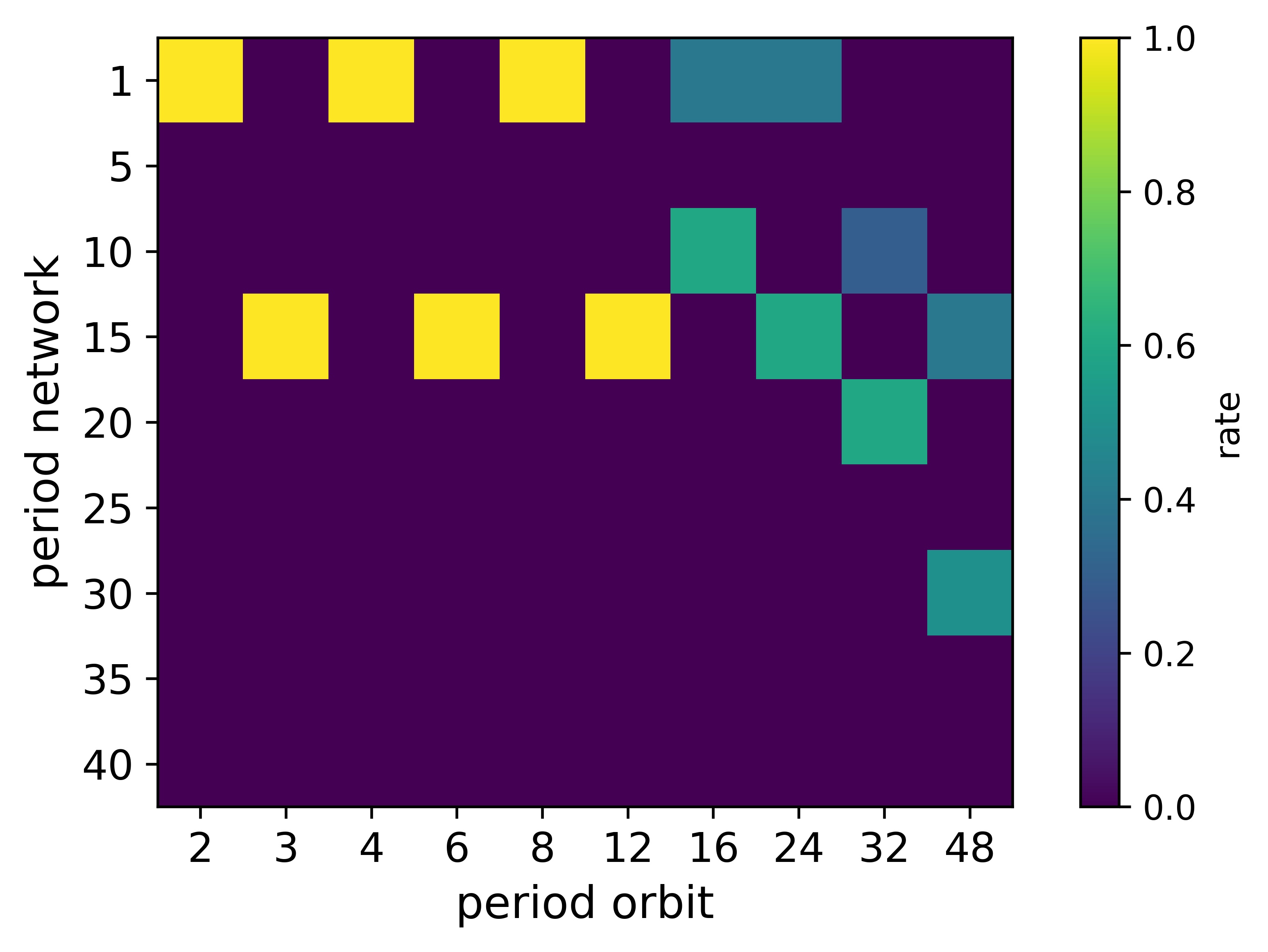}
         \caption{ 1\%  \label{fig:85raw}}
     \end{subfigure}
\caption{ Non-binary period matrices: The percentage indicates the size of the class with respect to all the models \label{fig:rawperiod}}     
\end{figure}

Recall that for low orbit periods we generally obtain a unique period of the network. Most high orbit periods map onto more than one network period. We normalize the period network over each period orbit by dividing by the total orbits of a specific period and then put the values in a matrix as before, see Figure \ref{fig:rawperiod}.

In Figure \ref{fig:rawperiod} we visualized the 4 most common classes. Figure \ref{fig:rawA} is the major class making up Figure \ref{fig:Aoverview}. We also that the max over each column of Figure \ref{fig:49raw} gives Figure \ref{fig:7max}. If we take the max over each column of Figure \ref{fig:52raw} or Figure \ref{fig:85raw} we do not get Figure \ref{fig:11max} or Figure \ref{fig:14max}.

\section{Performance models non-periodic sine-circle data \label{app:perform}}

\begin{figure}[h]
    \centering
    \includegraphics[width=7cm]{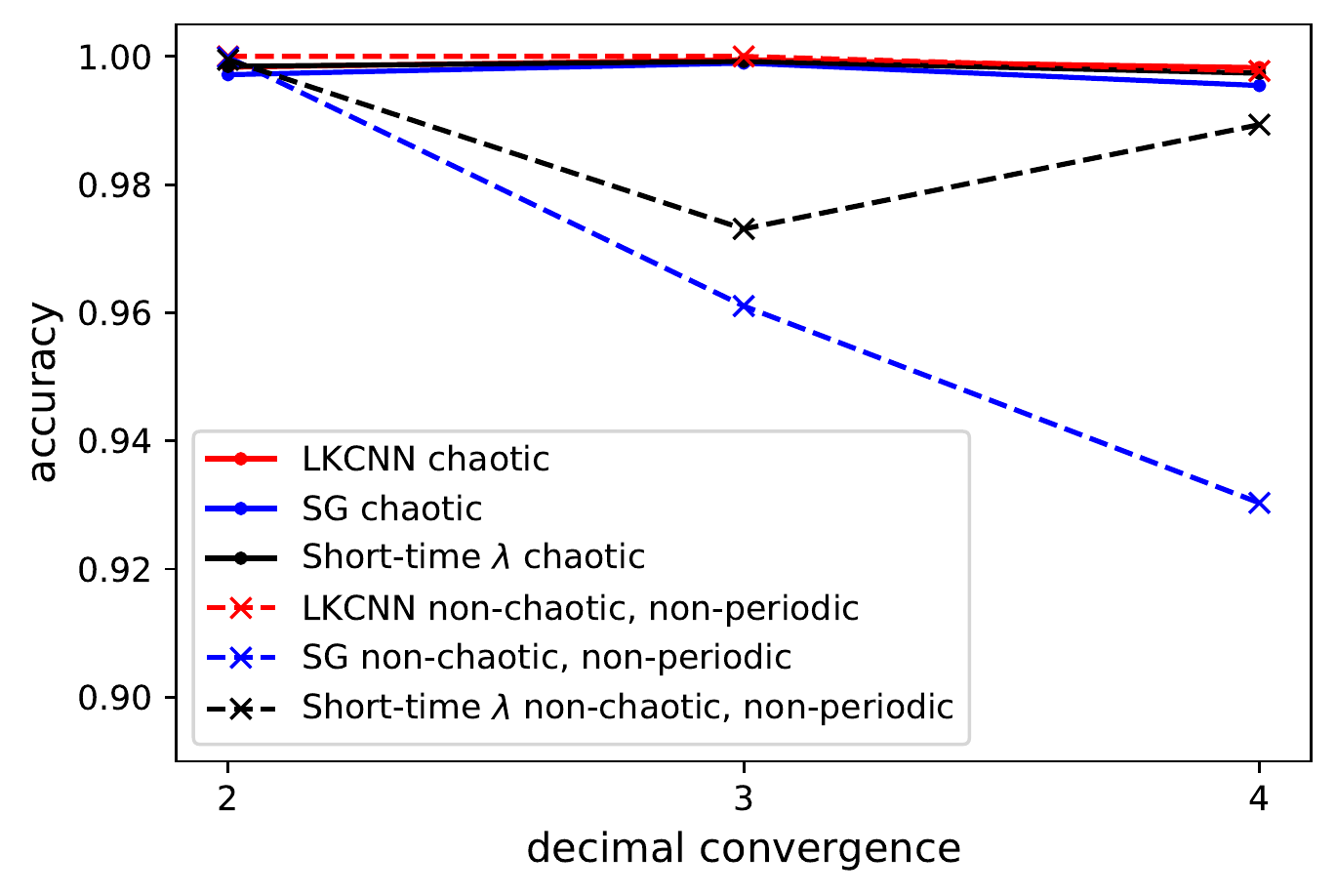}
    \caption{Performance models on non-periodic sine-circle data set as function of decimal convergence of Lyapunov exponent: We consider the performance of LKCNN, Savitsky-Golay reconstruction of the Lyapunov exponent (SG) and short-time $\lambda$ (averaged Lyapunov exponent over the input sequence). The $x$-axis corresponds to a subset where the first $k$-decimals of the Lyapunov exponent converge. Hence, the classification will also be determined by the first $k$-decimals. All models have near perfect performance when we only consider 2-decimal convergence or restrict to the chaotic subset. \label{fig:sg_perform}}
\end{figure}

In Table \ref{tab:diff_models} of Section \ref{sec:perform_LKCNN} we presented the performance of  LKCNN, Savitsky-Golay reconstruction of the Lyapunov exponent and short-time Lyapunov exponent on classifying chaos for non-periodic trajectories of the sine-circle map. For non-periodic trajectories the performance depends strongly on the decimal precision of the Lyapunov exponent. Here we consider sequences which have a Lyapunov exponent that converges in $k$-decimals and classify the sequence as chaotic using the first $k$-decimals. In Figure \ref{fig:sg_perform} we consider the performance as function of the convergence of the first $k$-decimals of the Lyapunov exponent.

The three models have near perfect performance when we only consider 2-decimal convergence or restrict to the chaotic subset. We note that the qualitative observations of Table \ref{tab:diff_models} also apply if we consider convergence in 3-decimals.  Evaluating the Lyapunov exponent over regular non-periodic sequences we typically observe fluctuations around zero which makes the classification task more difficult which is clearly reflected when we compare the performance for decimal precision greater than 2. Observe that performance of short-time $\lambda$ is non-monotone over the regular non-periodic subset. The sequences corresponding to $k$-decimal convergence are a subset of the sequences corresponding to $(k+1)$-decimal convergence. This property is insufficient to conclude anything about monotonicty in the performance results. These subsets do decrease in size with increasing $k$. For $k=5$ the resulting subset is so small that the results are not reliable from a generalization perspective.  

\section{Regular non-periodic logistic map trajectories}

If we train the LKCNN on the sine-circle data set from Section \ref{sec:data} we obtain poor performance on the regular non-periodic set, see Table \ref{tab:perform}. This begs the question why we obtain such good performance on the logistic set for the regular non-periodic subset.  It turns out that the majority of this subset is in a sense close to periodic trajectories of low period. We divide $x$ in consecutive chunks of length $k$ by defining 
\[
y_k^m(x) := (x_{1+k(m-1)} , \ldots, x_{1+km}).
\]
We consider the error given by minimizing over $k$ the average difference between consecutive  $y_k^m(x)$-chunks:
\begin{align*}
{\rm period \; error}(x) &:= \min_{k \leq \frac{n}{2}}{\frac{1}{\lfloor n/k \rfloor }\sum^{ \lfloor n/k \rfloor }_{m=1} \left| y_k^{m-1}(x)  - y_{k}^m(x)  \right|}, \\
K(x) &:= {\rm min} \; { \underset{k \leq \frac{n}{2}}{\rm argmin} {\frac{1}{\lfloor n/k \rfloor }\sum^{ \lfloor n/k \rfloor }_{m=1} \left| y_k^{m-1}(x)  - y_{k}^m(x)  \right|}}.
\end{align*}
The period error measures how close $x \in \mathbb{R}^n$ is to a periodic orbit.  More specifically, we have that if $\{\tilde{x}_i\}$ is a $k$-periodic sequence and $x=(\tilde{x}_1,\tilde{x}_2, \cdots, \tilde{x}_n )$ with $n\geq 2k$ then ${\rm period \; error}(x)=0$ and $K(x) = k$.  

Note that it would be computationally unfeasible to compare all chunks. However, for the computation we consider the minimum of the period error for the original and reversed sequence.

\begin{figure}[h]
    \begin{center}
        \includegraphics[width=7cm]{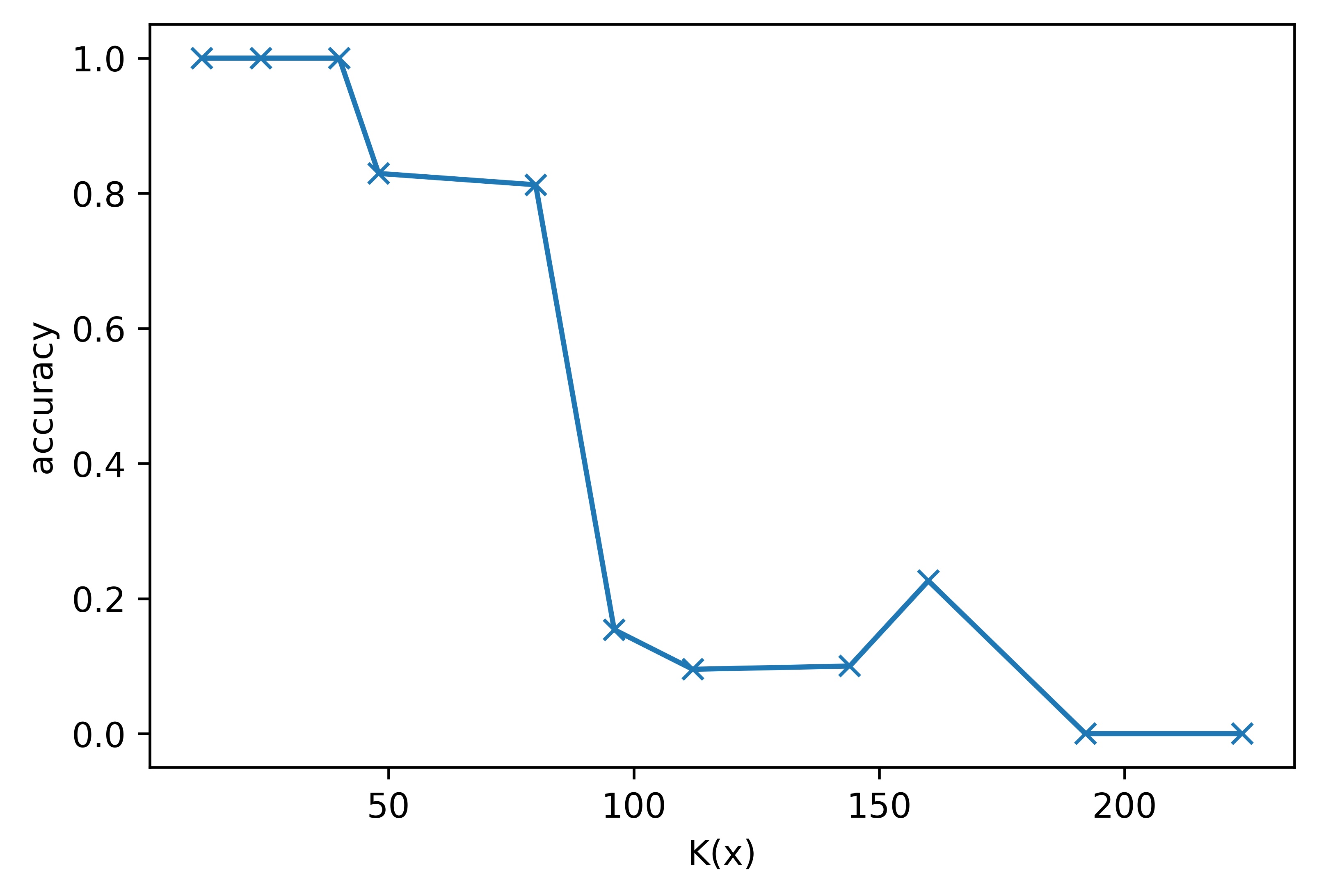}
    \end{center}
    \caption{Accuracy non-periodic trajectories logistic data: As the approximate period $K(x)$ is increased the accuracy decreases. This gives evidence that the network can only classify these trajectories correctly if they are close to trajectories of sufficiently small period. An exception here are periods which are a power of 2. These periods have been excluded in the graph since the network has accuracy close to 1 on these subsets. This is to be expected since the majority of periodic trajectories have period $2^k$} \label{fig:hist_nonper} 
\end{figure}

\begin{figure}[ht]
    \centering
    \includegraphics[width=7cm]{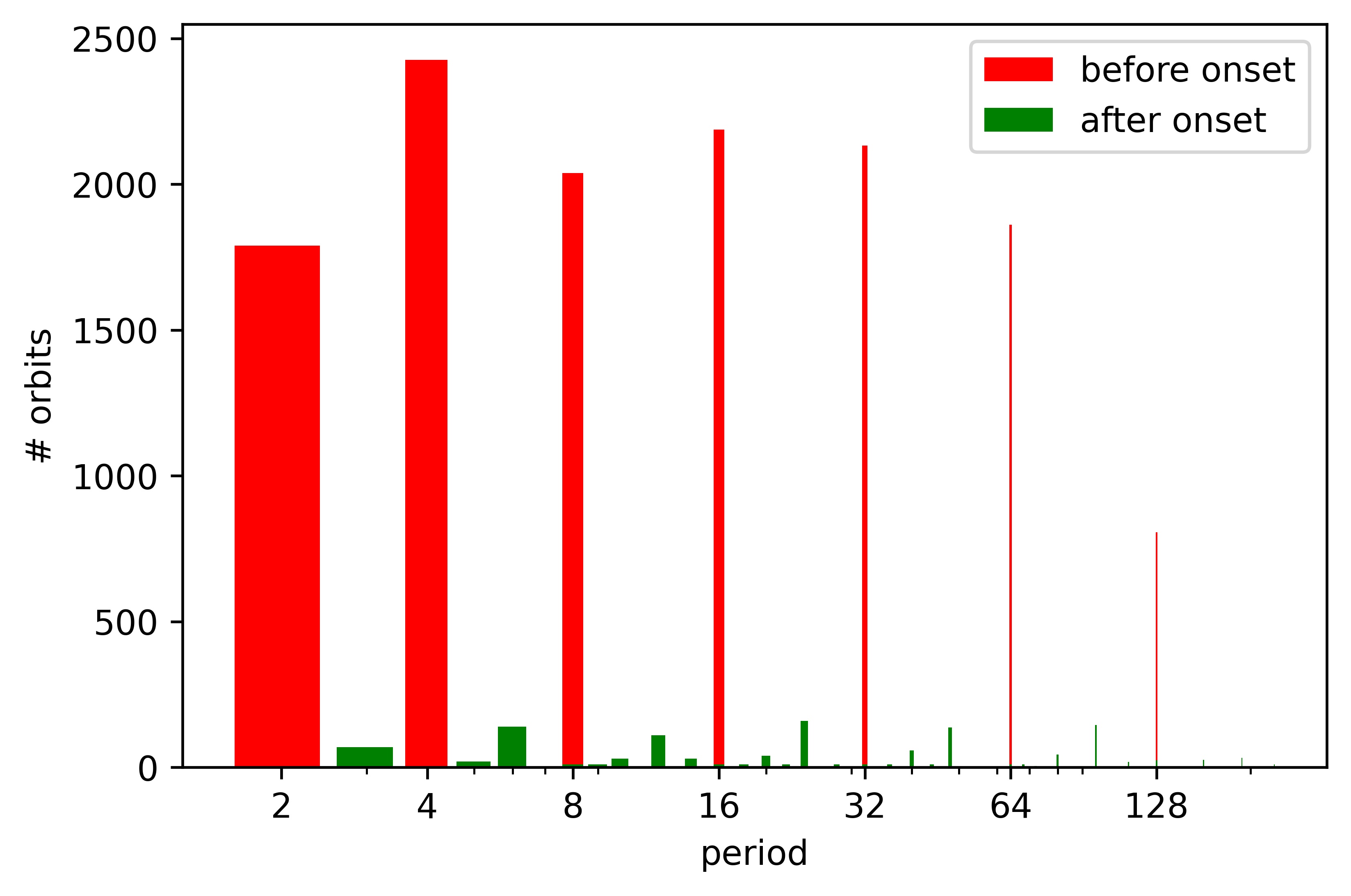}
    \caption{Log-spaced period distribution for logistic data set.}
    \label{fig:log-spaced2}
\end{figure}

For ${\rm period \; error} (x) < 10^{-5}$ we consider the accuracy of non-periodic orbits which are labeled as regular. If $K(x)=2^k$ for $k<8$ the network scores high accuracy since the majority of periodic orbits in the data have period $2^k$, see Figure \ref{fig:log-spaced2}.

In Figure \ref{fig:hist_nonper} we exclude $K(x)=2^k$. We observe that only for non-periodic trajectories with small $K(x)$ we obtain high-accuracy.

\bibliographystyle{alpha}

\bibliography{my_bib_abbreviated.bib}

\end{document}